\documentclass{article}
\usepackage[a4paper, margin=2.5cm]{geometry}
\usepackage[title]{appendix}

\usepackage{hyperref}

\usepackage{amsmath}
\usepackage{amssymb}

\usepackage{algorithm}    
\usepackage[noend]{algpseudocode}
  \algnewcommand\Yield{\State \textbf{yield}\ }
  \algnewcommand\Precondition{\State \textbf{Precondition:}\ }
  \algnewcommand\Input{\State \textbf{Input:}\ }
  
\usepackage{tikz}
  \usetikzlibrary{calc,positioning,backgrounds}
  \usetikzlibrary{
    decorations.pathmorphing,
    decorations.pathreplacing,
    decorations.shapes
   } 
  \usetikzlibrary{arrows.meta, chains, fit, quotes, shapes.symbols, shapes.misc, shapes.geometric}
  \usetikzlibrary {shadows}
  \usetikzlibrary{patterns, patterns.meta}

\usepackage{threeparttable}
\usepackage{booktabs}
\usepackage{subcaption}
\usepackage{cleveref}

\newcommand{\action}[1]{\tilde{#1}}
\newcommand{\tnow}{t^{\operatorname{now}}}

\usepackage{natbib}
\setcitestyle{authoryear,open={(},close={)}}

\tikzset{
  fc-state/.style = {
    rectangle,
    draw= black, fill= white,
    minimum width= 3cm, minimum height= 0.8cm,
    align=center
  },  
  fc-quest/.style = {
    chamfered rectangle,
    chamfered rectangle xsep=3cm,
    draw= black, fill= white,
    align=center
  },
  fc-startstop/.style = {
    rounded rectangle,
    rounded rectangle arc length=180,
    draw= black, fill= white,
    align=center,
    minimum width= 2cm, minimum height= 0.6cm,
    scale=0.8
  },
  fc-arrow/.style = {
    draw,
    ->, -latex,
    thick 
  },
  fc-arrowlabel/.style = {
    rectangle,
    fill=white,
    thin,
    align= center,
    scale= 0.8
  },
  fc-fit/.style={
    rectangle,
    draw=black!50!white, fill=black!2!white
  },
  dp-net/.style={
    regular polygon,
    regular polygon sides=6,
    minimum size= 3cm,
  },
  dp-veh/.style={
    draw=black,
    fill=white,
    rectangle,
    minimum height= 0.5cm,
    minimum width=  0.5cm
  },
  dp-loc/.style={
    draw=black,
    fill=white,
    circle,
    minimum height= 0.125cm,
    minimum width=  0.125cm
  },
  dp-locin/.style={
    draw=black,
    fill=black,
    circle,
    minimum height= 0.125cm,
    minimum width=  0.125cm
  },
  dp-loclab/.style={
    draw=none,
    fill=none,
    scale=0.8
  },
  dp-arrow/.style={
    draw,
    ->, -latex,
    thick
  },
  dp-orderarrow/.style={
    dp-arrow,
    dashed
  },
  dp-arrowlabel/.style={
    draw=none,
    circle,
    fill=white,
    thin,
    align= center,
    scale= 0.8,
    midway
  }
}

\newcommand{\sState}{\mathcal{S}}

\newcommand{\sVehStatus}{\Phi}
\newcommand{\sOrdStatus}{\Psi}
\newcommand{\sLoad}{\mathcal{C}}

\title{A general modeling and simulation framework for dynamic vehicle routing}

\author{Mark\'o Horv\'ath\thanks{HUN-REN Institute for Computer Science and Control, Budapest, Hungary; marko.horvath@sztaki.hu; corresponding author}
\and
Tímea Tam\'asi\thanks{Department of Operations Research, Institute of Mathematics, ELTE Eötvös Loránd University, Budapest, Hungary and HUN-REN Institute for Computer Science and Control, Budapest, Hungary; tamasitimea@student.elte.hu}
}

\begin{document}

\maketitle

\begin{abstract}
In dynamic vehicle routing problems (DVRPs), some part of the information is revealed or changed on the fly, and the decision maker has the opportunity to re-plan the vehicle routes during their execution, reflecting on the changes.
Accordingly, the solution to a DVRP is a flexible policy rather than a set of fixed routes.
A policy is basically a problem-specific algorithm that is invoked at various decision points in the planning horizon and returns a decision according to the current state.
Since DVRPs involve dynamic decision making, a simulator is an essential tool for dynamically testing and evaluating the policies.
Despite this, there are few tools available that are specifically designed for this purpose.
To fill this gap, we have developed a simulation framework that is suitable for a wide range of dynamic vehicle routing problems and allows to dynamically test different policies for the given problem.
In this paper, we present the background of this simulation tool, for which we proposed a general modeling framework suitable for formalizing DVRPs independently of simulation purposes.
Our open source simulation tool is already available, easy to use, and easily customizable, making it a useful tool for the research community.\\

\end{abstract}

\section{Introduction}

A vehicle routing problem is \emph{dynamic}, if some part of the information is revealed or changed on the fly, and the decision maker (the service provider) has the opportunity to re-plan the vehicle routes during their execution, reflecting on the changes.
Dynamic vehicle routing problems (DVRPs) have received a lot of attention in the past decades, which is certified by a series of recent review papers, e.g., \citet{berbeglia2010dynamic, pillac2013review, bekta2014chapter, psaraftis2016dynamic, ritzinger2016survey, rios2021recent, soeffker2022stochastic, zhang2023dynamic, mardevsic2023review}.
This growing interest is due to the wide range of real-world applications and the fact that today's technology enables real-time decision making.

\begin{figure}
\tikzset{
  process/.style = {
    signal, signal to=east, draw, right, fill=white
  },
  eventnode/.style = {
    circle, draw=black, fill=white, inner sep=0, minimum size=0.2cm
  },
  decisionnode/.style = {
    eventnode, fill=black
  },
  state/.style = {
    rectangle, draw, fill=white, inner sep=0, minimum size= 0.2cm
  },
  pdstate/.style = {
    state, fill=black
  }
}
\centering
\begin{subfigure}[t]{0.45\textwidth}
\centering
\begin{tikzpicture}

\node[process, minimum width= 5cm, minimum height= 0.5cm, above right= 0.0cm and 0cm] (sdp) at (0,0) {};
\node[scale=0.8,font=\it,above= 0.25cm of sdp] (sdplabel) {decision process};

\node[decisionnode]         (n1) at (1.1,0.25) {};
\node[eventnode,draw=none]  (n2) at (1.7,0.25) {};
\node[decisionnode]         (n3) at (3.6,0.25) {};
\node[eventnode,draw=none]  (n4) at (4.1,0.25) {};

\node[state,below=0.5 of n1] (s1) {};
\node[pdstate,below=0.5 of n2] (s2) {};
\node[state,below=0.5 of n3] (s3) {};
\node[pdstate,below=0.5 of n4] (s4) {};
\node[state,left=0.85 of s1] (s0) {};

\foreach \i in {1,...,4}{
  \draw[fc-arrow] (s\the\numexpr\i-1\relax) -- (s\i);
}

\foreach \i in {1,3}{
  \draw[fc-arrow] (n\i) -- (s\i);
}

\node[process, minimum width= 5cm, minimum height= 0.5cm] (dm) at (0,-2.5) {};
\node[scale=0.8,font=\it] (dmlabel) at (dm.center) {decision making};

\draw[fc-arrow] (s1) -- node[scale=0.8,fc-arrowlabel,pos=0.33,draw] (ls1) {state} (s1 |- dm.north);
\draw[fc-arrow] (s2 |- dm.north) -- node[scale=0.8,fc-arrowlabel,pos=0.33,draw] {action} (s2);
\draw[fc-arrow] (s3) -- node[scale=0.8,fc-arrowlabel,pos=0.33,draw] (ls2) {state} (s3 |- dm.north);
\draw[fc-arrow] (s4 |- dm.north) -- node[scale=0.8,fc-arrowlabel,pos=0.33,draw] {action} (s4);

\scoped[on background layer]{
\node[fc-fit,fit=(sdp)(sdplabel)(s0)]{};
}
\end{tikzpicture}
\caption{Sequential decision process.}
\label{fig:sdp:sdp}
\end{subfigure}
\begin{subfigure}[t]{0.45\textwidth}
\centering
\begin{tikzpicture}
\node[process, minimum width= 5cm, minimum height= 0.5cm, above right= 0.0cm and 0cm] (sdp) at (0,0) {};
\node[scale=0.8,font=\it,above= 0.25cm of sdp] (sdplabel) {decision process};

\node[eventnode]    (n1) at (0.5,0.25) {};
\node[decisionnode] (n2) at (1.1,0.25) {};
\node[eventnode,draw=none]    (n3) at (1.7,0.25) {};
\node[eventnode]    (n4) at (2.4,0.25) {};
\node[eventnode]    (n5) at (3.0,0.25) {};
\node[decisionnode] (n6) at (3.6,0.25) {};
\node[eventnode,draw=none]    (n7) at (4.1,0.25) {};
\node[eventnode]    (n8) at (4.6,0.25) {};

\node[state,below=0.5 of n1] (s1) {};
\node[state,below=0.5 of n2] (s2) {};
\node[pdstate,below=0.5 of n3] (s3) {};
\node[state,below=0.5 of n4] (s4) {};
\node[state,below=0.5 of n5] (s5) {};
\node[state,below=0.5 of n6] (s6) {};
\node[pdstate,below=0.5 of n7] (s7) {};
\node[state,below=0.5 of n8] (s8) {};
\node[state,left=0.25 of s1] (s0) {};

\foreach \i in {1,...,8}{
  \draw[fc-arrow] (s\the\numexpr\i-1\relax) -- (s\i);
}

\foreach \i in {1,2,4,5,6,8}{
  \draw[fc-arrow] (n\i) -- (s\i);
}

\node[process, minimum width= 5cm, minimum height= 0.5cm] (dm) at (0,-2.5) {};
\node[scale=0.8,font=\it] (dmlabel) at (dm.center) {decision making};

\draw[fc-arrow] (s2) -- node[scale=0.8,fc-arrowlabel,pos=0.33,draw] (ls1) {state} (s2 |- dm.north);
\draw[fc-arrow] (s3 |- dm.north) -- node[scale=0.8,fc-arrowlabel,pos=0.33,draw] {action} (s3);
\draw[fc-arrow] (s6) -- node[scale=0.8,fc-arrowlabel,pos=0.33,draw] (ls2) {state} (s6 |- dm.north);
\draw[fc-arrow] (s7 |- dm.north) -- node[scale=0.8,fc-arrowlabel,pos=0.33,draw] {action} (s7);

\scoped[on background layer]{
\node[fc-fit,fit=(sdp)(sdplabel)(s0)]{};
}
\end{tikzpicture}
\caption{Discrete-event based decision process.}
\label{fig:sdp:dep}
\end{subfigure}
\caption{
Differences between the sequential and the discrete-event based decision process.
Circles refer to distinct events (e.g., order requests, vehicle arrival).
Black circles refer to decision points.
Squares refer to states.
Black squares refer to post-decision states.
}
\end{figure}
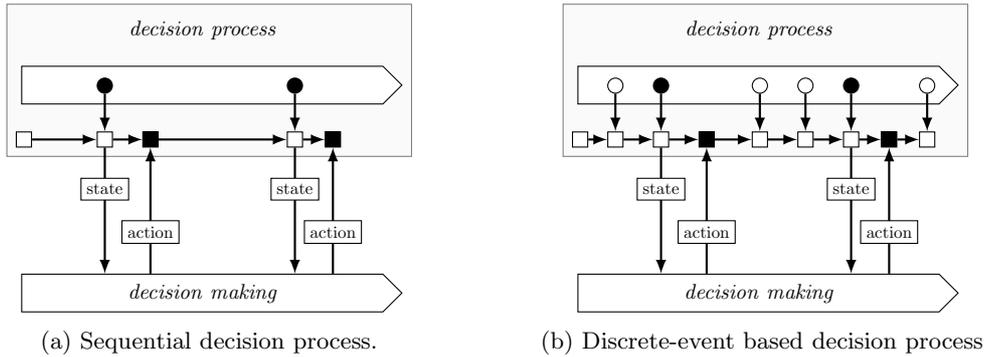

Nowadays, DVRPs are usually modeled using the so-called \emph{sequential decision process} (e.g., \citet{ulmer2020modeling, soeffker2022stochastic}).
Briefly stated, the decision process transitions from decision point to decision point, where the decision maker is provided with the current state (i.e., all the available information) and has the opportunity to make a decision (e.g., update the vehicle routes), or in other words, to choose an action, see \Cref{fig:sdp:sdp}.
Accordingly, a solution to the dynamic problem is a \emph{policy}, which is a function that assigns an action to every state.

Apart from survey articles, in the majority of the papers dealing with DVRPs, the authors propose policies for the problem at hand, and perform computational experiments to evaluate them, e.g., to compare them with state-of-the-art or baseline policies.
In addition to doing the obviously necessary implementation of their policy, they need some kind of simulator for dynamic evaluation.
In this paper, we focus on this dynamic evaluation, and we approach the DVRPs from the simulation point of view.
Even more emphasized, our focus is not on the solution approaches for a particular DVRP, but on the modeling of general problems and on the dynamic testing of arbitrary solution methods.

According to our primary goal, we have implemented a simulation framework that is suitable for a wide range of dynamic vehicle routing problems and allows to dynamically test different solution approaches for the modeled problem.
This article, however, is much more than technical documentation, as we also propose a general modeling framework suitable for formalizing DVRPs independently of simulation purposes.

\subsection{Motivation}

The simulation of the decision process is essential for the dynamic evaluation of solution approaches to dynamic vehicle routing problems.
Despite this, there are few tools available that are \emph{specifically} designed for this purpose.

In simpler cases, it is very easy to implement the sequential decision process, since the transition between the states is straightforward.
However, in many other cases (especially when inter-route constraints make the problem difficult), it is necessary to run a more complex simulation to move the decision process from decision point to decision point.
Although general-purpose simulation tools exist (e.g., AnyLogic, SimPy), they require the user to build the entire dynamic vehicle routing framework from scratch.
Several publicly available simulators have been created using these tools, but they are only suitable for a specific problem (see e.g., \citet{hao2022introduction}).
Transportation simulation software packages (e.g., Eclipse SUMO, MATSim, PTV Vissim, Transims) could potentially support dynamic testing, but most of these tools focus primarily on microscopic traffic simulation (including elements such as traffic lights and pedestrian interactions), a level of detail that is rarely considered for research in our scope.
We would like to highlight the work of \citet{maciejewski2016dynamic, maciejewski2017towards}, where the authors developed a DVRP extension for MATSim.
This extension allows the modeling of a wide variety of DVRPs and the plugging of different algorithms, therefore this tool is indeed suitable for dynamical testing.
However, modeling and customization requires familiarity with Java and the relatively complex architecture of MATSim, including a batch of scenario files.
Our understanding is that the implementation of the decision making algorithm is also tied to Java.

Based on the above, it is a reasonable goal to develop a standalone simulation tool for DVRPs according to the following criteria.
(i) The simulation tool should be based on a generic modeling framework in which the problems can be clearly formulated, thus ensuring the reconstruction of the research.
(ii) The framework should be able to model a wide range of DVRPs, that is, the problem aspects and side constraints often occur in the literature should be included by default.
(iii) The simulation tool should be easy to use, so it should be much easier to model a problem in it than to implement an entire decision process from scratch.
(iv) The simulation tool should be easily customizable and adaptable to individual needs.
(v) The implementation of the decision making algorithm should not be tied to a specific programming language, but the simulator should allow communication with it.

\subsection{Main contributions}

Our main goal was to develop a general simulation tool for dynamic vehicle routing.
To achieve this, we conducted an extensive literature review and developed a general modeling and simulation framework.
Our main contributions are the following.

\paragraph{Literature review on DVRPs}
We studied the literature on dynamic vehicle routing to identify those problem aspects and side constraints that are common and should therefore be considered in the development of the framework.
For details, see \Cref{sec:lit_rev}.

\paragraph{Modeling and simulation framework for DVRPs}
We developed a general modeling and simulation framework for dynamic vehicle routing.
The framework is suitable for modeling a wide range of DVRPs, primarily pickup-and-delivery problems, but it is easily adaptable to other problems as well.
Our \emph{discrete-event based decision process} is a combination of the discrete-event based simulation and the sequential decision process, the latter of which is widely used to formalize DVRPs.
For the modeling, we borrowed the route-based representation of \citet{ulmer2020modeling}, but we propose a more detailed model suitable for simulation purposes, see \Cref{fig:sdp:dep}.
We also standardized and formalized some common aspects of decision making, such as postponing decisions and delaying the departure of vehicles.
For details, see \Cref{sec:basics,sec:mod}.

\paragraph{Open source simulation tool for DVRPs}
According to our primary goal, we created an implementation of our simulation framework.
The source code of our Python package, called \texttt{dvrpsim}, is available online.
Dynamic vehicle routing problems can be easily modeled, and the simulator is easily customizable, making it a useful tool for other researchers to dynamically test and evaluate their algorithms for a particular problem.
To the best of our knowledge, this is the first simulation tool designed specifically for this purpose.
For details, see \Cref{sec:simulation}.

\section{Dynamic vehicle routing}\label{sec:lit_rev}

In this section, we provide a brief introduction to dynamic vehicle routing.
We also summarize our literature review on dynamic vehicle routing problems.
We compiled the reviewed papers in \Cref{tab:lit:vrpdsr,tab:lit:dpdp,tab:lit:sddp}.
The goal of the review was to identify those problem aspects and side constraints that often occur in the literature, therefore, they should be taken into account when developing a general modeling and simulation framework.
As the focus is on modeling and simulation, the literature review does not cover problem aspects such as logistic context, objective functions, solution approaches, etc.
For such an overview, we refer to the excellent review by \citet{zhang2023dynamic}.

\subsection{Dynamic vehicle routing problems}\label{sec:dvr:dvrps}

Briefly stated, the well-known \emph{(static) vehicle routing problem} (VRP) aims to determine an optimal set of routes to be performed by a fleet of vehicles to fulfill order requests at different locations within a planning horizon.
The problem was introduced more than 60 years ago by \citet{dantzig1959truck}, then generalized by \citet{clarke1964scheduling}, and many variations have appeared since then (e.g., \citet{toth2002vehicle, eksioglu2009vehicle, braekers2016vehicle, zhang2022review}).

According to \citet{psaraftis1980dynamic}, a vehicle routing problem is characterized as \emph{dynamic}, if the input of the problem is received and updated concurrently with the determination of the routes.
The vehicle routes can be redefined in an ongoing fashion.
This class of problems is often referred to as \emph{online} or \emph{real-time}.
Using the taxonomy of \citet{pillac2013review}, a dynamic problem is \emph{stochastic}, if there is some exploitable stochastic knowledge about the dynamically revealed information, and \emph{deterministic} otherwise.
Thus, \emph{stochastic dynamic vehicle routing problems} (SDVRPs) are also within the scope of our paper.

In a recent survey, \citet{zhang2023dynamic} considered three DVRP subcategories by distinguishing three types of order requests.
(i)~A \emph{pickup and delivery request} consists of a pair of locations, and the serving vehicle must visit the pickup location before going to the delivery location.
\Cref{tab:lit:dpdp} summarizes the papers we have reviewed on the associated \emph{dynamic pickup-and-delivery problems} (DPDPs).
(ii)~\emph{Delivery requests} are special pickup and delivery requests because their pickup location refers to a depot.
See \Cref{tab:lit:sddp} for our summary on the related \emph{same-day delivery problems} (SDDPs).
(iii)~A \emph{service request} is associated with only a single location, so the assigned vehicle does not have to visit a specific pickup location (e.g., the depot) before serving the request.
See \Cref{tab:lit:vrpdsr} for our overview on \emph{vehicle routing problems with dynamic service requests} (VRPDSRs).

\paragraph{Problems in our scope}

In this paper, we focus on the three DVRP subcategories considered by \citet{zhang2023dynamic}.
We present our modeling framework primarily for DPDPs (including SDDPs) as we assume that each request has a designated origin and a designated destination, however, with a slight modification the framework is also adaptable to DVRPs with service requests.

Note that \citet{zhang2023dynamic} identified another DVRP variant in addition to the previous ones, called the \emph{dynamic multi-period VRP} (DMPVRP), which is characterized by multiple planning periods.
In this paper, we do not consider these problems.
We also do not consider those problems, where the transportation consists of multiple stages, such as \emph{multi-echelon vehicle routing} or \emph{vehicle routing with transshipment}.
For a review on these problems, see e.g., \citet{sluijk2023two, nielsen2024systematic}.

\subsection{Sequential decision process}

Nowadays, the state-of-the-art approach to modeling DVRPs is the \emph{sequential (or Markov) decision process}.
For a thorough introduction, see \citep{ulmer2020modeling, soeffker2022stochastic}.
Briefly stated, at certain time points in the planning horizon, called \emph{decision points}, the decision maker has the opportunity to re-plan the vehicle routes, reflecting on the newly revealed information, see \Cref{fig:sdp:sdp}.
These decision points may be predetermined (e.g., they occur at given intervals), or they can be imposed by certain events (e.g., requesting an order).
The sequential decision process steps from decision point to decision point, called \emph{transition}.
At a decision point, the decision maker is provided with the current \emph{state}, which describes all the information available to make a decision.
The result of the decision making is an \emph{action} that includes, for example, the updated vehicle routes.

Note in advance that our discrete-event based decision process differs from the sequential decision process in that it explicitly considers events between decision points, see \Cref{fig:sdp:dep}.
Besides the fact that this approach makes it easier to formalize the dynamic problem in some cases, this level of detail allows us to construct a general, easily customizable simulation framework.

\subsection{Problem aspects and side constraints}\label{sec:bas:problem}

Now, we present the main aspects and side constraints of dynamic vehicle routing problems that were considered when building our framework.
We group these aspects by locations (\Cref{sec:bas:prob:locs}), orders (\Cref{sec:bas:prob:ords}), and vehicles (\Cref{sec:bas:prob:vehs}), but there may be some overlap between the groups.

\subsubsection{Locations}\label{sec:bas:prob:locs}

Location is a collective term for the places that vehicles may visit, such as depots, customers, restaurants, factories, etc., depending on the problem at hand.

\paragraph{Operating network}

At this level of logistics planning, vehicles operate on networks.
That is, the movement of vehicles is not detailed; they are either at a location (residing at a network node) or on the way (traveling along a network edge).
In the latter case, the exact positions of the vehicles are unknown, but their arrival can be calculated from the travel time.
In certain cases, vehicle movements are simulated within a real-world road network, such that road crossings also refer to locations (e.g., \citet{ferrucci2014realtime, ferrucci2015general, ferrucci2016proactive}).
Vehicles, especially if they are different types (e.g. drones and trucks), can operate on different networks  (e.g., \citet{ulmer2018sameday}).

\paragraph{Travel times}

Travel times between locations can be vehicle-dependent, for example, if vehicles have different speeds, and especially if the vehicles operate on different networks (e.g., \citet{ulmer2018sameday}).
Travel times can also be time-dependent (e.g., \citet{haghani2005dynamic}) or even stochastic (e.g., \citet{schilde2014integrating}).

\paragraph{Docking restrictions}

Locations often have limited space for loading or unloading, and sometimes the loading crew creates a bottleneck.
Because of these inter-route constraints, vehicles may make each other wait.
For example, \citet{hao2022introduction} proposed a problem, where each factory has a limited number of docking ports, so if a vehicle arrives and there is no port available, the vehicle must wait until a port becomes available.

\subsubsection{Orders}\label{sec:bas:prob:ords}

Orders are transportation or service requests.
The object of transportation can be a variety of products, food (e.g., meal delivery problem), other vehicles (e.g., bike sharing rebalancing problem), or even people (e.g., dial-a-ride problem).
Transportation requests typically have an origin (i.e., pickup location) and a destination (i.e., delivery location).
In many cases, the terms "order" and "customer" are used interchangeably.

\paragraph{Service times}

Various service times may arise when orders are picked up or delivered.
Loading and unloading itself may take some time and may even depend on the quantity of orders (e.g., \citet{hao2022introduction}).
These times can also be location-dependent (e.g., \citet{ulmer2019preemptive}) or vehicle-dependent (e.g., \citet{ulmer2018sameday}).
Additional order-independent service times, such as parking or docking, may also occur (e.g., \citet{hao2022introduction}).
The above service times can even be stochastic (e.g., \citet{goel2019vehicle}), but in many cases they are simply neglected or incorporated in the travel times.

\paragraph{Service time windows}

Orders often have \emph{service time windows} for their pickup and/or delivery.
Such a time window specifies an \emph{earliest} and a \emph{latest service start time} for the order.
Earliest service start times are typically hard constraints, meaning that if a vehicle arrives early at a location, it has to wait until the time window opens, but \citet{schilde2014integrating}, for example, allowed early arrivals.
In contrast, latest service start times are often soft constraints, that is, the service can start after the latest required time, however, the tardiness may incur additional costs (e.g., \citet{ulmer2021restaurant}).
In some rare cases, customers have multiple time windows in the planning horizon (e.g., \citet{de2015constrained, de2015variable}).
Service time windows can be stochastic.
For example, in the problem proposed by \citet{srour2018strategies}, customers first preannounce their request with an estimated time window for pickup, which can be changed when the customer confirms the request.

\paragraph{Order cancellation}

In some cases, customers can cancel their requests (e.g., \citet{lin2014decision, los2020value}).
Cancellation is allowed only if the service of the corresponding order has not yet started.
After the notification, the decision maker must remove the canceled orders from the vehicle routes.
Cancellation is permanent, and canceled orders are no longer dealt with in the given planning horizon.

\subsubsection{Vehicles}\label{sec:bas:prob:vehs}

Vehicle is a collective term for the equipment or people that perform the transportation, such as trucks, drones, drivers, couriers, etc., depending on the problem at hand.

\paragraph{Vehicle fleet}

The fleet of vehicles can be either \emph{homogeneous} or \emph{heterogeneous}.
In the latter case, vehicles may differ not only in their basic parameters, but also in their operations.
For example, \citet{ulmer2018sameday} considered a problem with heterogeneous fleets of drones and trucks that differ not only in their availability, capacity, and travel speed, but also in their requirement for charging and the network on which they operate.

\paragraph{Vehicle capacity}

A vehicle is either \emph{capacitated} or \emph{uncapacitated}.
In the former case, the total size or quantity of orders loaded on the vehicle must never exceed the capacity of the vehicle.
In dial-a-ride or taxi-routing problems, the capacity of the vehicles is the number of non-driver seats, however, in some cases no shared rides are allowed, that is, a vehicle can only carry one passenger (or one passenger group) at a time (e.g., \citet{hyland2018dynamic}).
The uncapacitated case is common with those problems where the packages are relatively small and therefore the trunk of the transporting vehicle is not a limiting factor.

\paragraph{Loading rule}

Vehicles can be subject to \emph{loading rules}.
For example, in \citep{hao2022introduction}, unloading must follow the last-in-first-out (LIFO) rule, i.e., the last loaded order must be unloaded first.

\paragraph{Vehicle availability}

Vehicles can also have time windows, representing the working shifts of the drivers (e.g., \citet{de2015variable, steever2019dynamic}).
Sometimes, a time window~$[0,L]$ is associated with the depot, also called the \emph{depot deadline}, which gives a latest return time~($L$) for the vehicles (e.g., \citet{cote2023branch}).


\subsection{Aspects in decision making}\label{sec:bas:decision}

Several questions may arise when making decisions.
When or how often is it necessary to re-optimize (\Cref{sec:int:dec:points})?
Can and should an order be rejected (\Cref{sec:int:dec:rejection})?
Should all decisions be taken as soon as possible, or can certain decisions be postponed (\Cref{sec:int:dec:opp})?
Should vehicles be sent on their way immediately or is it worth waiting (\Cref{sec:int:dec:delay})?
Can en route vehicles be diverted or should their destination not be changed (\Cref{sec:int:dec:enroute})?
Can a request be served by multiple routes (\Cref{sec:int:dec:split})?

\subsubsection{Decision points}\label{sec:int:dec:points}

In the case of DVRPs, the decision maker must decide when to process the new dynamic information and update the routes of the vehicles.
Most of the articles use three different approaches, namely the decision maker makes a decision either periodically, when a new order request arrives, or when a vehicle arrives at a location, however, there are several other possibilities, and the various approaches can also be combined.

\paragraph{Periodic decision points}

In many applications, the planning horizon is divided into predetermined decision epochs, typically of equal length ($\Delta$), i.e., decision points occur periodically.
For example, \citet{zolfagharinia2014benefit} re-planned truck routes twice a day ($\Delta = 12\,h$).
In the framework proposed by \citet{hao2022introduction} for a dynamic pickup-and-delivery problem, information is updated every 10 minutes ($\Delta = 10\,min$).
\citet{bertsimas2019online} re-optimized taxi routes even more frequently ($\Delta = 30\,s$).

\paragraph{Decision point on order request}

The most common case is that decisions are made when new orders are requested.
\citet{ninikas2014reoptimization} also considered a policy where, instead of imposing decision points on every order request, re-optimization would occur after a pre-defined number of requests.

\paragraph{Decision point on vehicle arrival}

Often, a decision point is imposed when a vehicle arrives at a location.
In some cases, complete order information is not available until arrival, so routes may need to be re-planned prior to the start of service (e.g., \citet{goodson2016restockingbased}).
For some same-day delivery problems, the planned vehicle routes are fixed, so re-optimization occurs only when a vehicle returns to the depot (e.g., \citet{dayarian2020sameday}).
In fact, most of the cases decision making is required after the service is finished, but since service times are neglected, it coincides with the arrival.
In many approaches, the planned route of a vehicle consists only of the next location to visit, so it is necessary to re-plan the route after the service is finished (e.g., \citet{ulmer2018budgeting, ulmer2019offlineonline}).

\paragraph{Self-imposed decision points}

In some cases, certain decisions can be postponed, which often involves the introduction of self-imposed decision points.
That is, if no other event imposes a decision point by a certain time point, then reaching that time will impose one to reconsider the decision.
For example, \citet{zhang2018dynamic} considered an orienteering problem in which a traveler must join a waiting queue upon arrival at a location.
If the traveler joins the queue, the next decision point is imposed when the size of the queue decreases or a predetermined maximum waiting time elapses, whichever occurs first.
\citet{ulmer2021restaurant} investigated a restaurant meal delivery problem, where the assignment of an order to a driver, once made, cannot be altered.
Thus, the authors proposed a policy, where the assignment of some non-urgent orders is postponed for a given unit of time, and if no new orders are requested during this period, the expiration of the postponement imposes a decision point.
In certain cases, delaying the departure of the vehicles can also cause self-imposed decision points, see later in \Cref{sec:int:dec:delay}.

\subsubsection{Order rejection}\label{sec:int:dec:rejection}

In many applications, the decision maker can reject orders, if they are unable or unwilling to fulfill them.
The rejection is permanent, and rejected orders are no longer dealt with in the given planning horizon.
In practice, rejected orders may be outsourced to a third party or moved to another planning horizon.
In the problem proposed by \citet{ehmke2014customer}, the decision maker allows the customer to request an alternative order with a different time window, if the original order is rejected.

\subsubsection{Decision postponement}\label{sec:int:dec:opp}

As we touched on in \Cref{sec:int:dec:points}, certain decisions can be postponed in some cases.
In our interpretation, decision postponement means that certain non-changeable decisions are not made at the current decision point, but are postponed to a later one.
For example, if order rejection is allowed, the acceptance/rejection is permanent, therefore some authors do not want to make the decision at the first possible decision point (e.g., \citet{zhang2018dynamic, voccia2019sameday}).
Sometimes, the assignment of orders to vehicles, once made, cannot be altered, so the decision on this assignment is postponed (e.g., \citet{ulmer2021restaurant}).
Note that the case where the order requests do not impose decision points, and the orders are accepted or rejected at the first decision point after their request, is not considered as decision postponement.

\subsubsection{Delaying the departure}\label{sec:int:dec:delay}

In addition to assigning routes to vehicles, it is also important to decide when to send vehicles on their way, since waiting for possible future orders could be beneficial.
The two basic waiting strategies, the \emph{drive-first} and the \emph{wait-first}, require a vehicle to departure from its current location at the earliest possible time and at the latest possible time, respectively, but several other waiting strategies have been applied to delay the departure of the vehicles (e.g., \citet{mitrovic2004waiting, branke2005waiting, ichoua2006exploiting}).

As mentioned in \Cref{sec:int:dec:points}, delaying the departure may involve the use of self-imposed decision points.
For example, \citet{voccia2019sameday} considered a same-day delivery problem, where the depot-to-depot tours cannot be modified during their execution.
In their policy, the authors did not start the vehicles immediately after determining their routes, but postponed them for a certain period of time.
A decision point was implied at the end of the waiting period, unless another event triggered one in the meantime.

\subsubsection{Diversion from the planned route}\label{sec:int:dec:enroute}

Due to the dynamic nature of the problem, the decision maker may modify the vehicle routes during execution.
Although the majority of papers consider decision making to be instantaneous, in practice it may cover longer periods of time during which the state of the system may change so much (e.g., some vehicles may have already departed) that the decision is no longer feasible with respect to this new state.
Therefore, it may be advisable to fix the first parts of the routes, i.e. to make them non-changeable.

In most SDDPs, once the vehicle leaves the depot, its entire route is fixed until it returns to the depot.
In some other cases, however, a \emph{preemptive depot return} is allowed, that is, the delivery vehicle can return to the depot before delivering all the orders it is currently carrying (e.g., \citet{ulmer2019preemptive, cote2023branch}).

In general, the next location of a vehicle is fixed.
This is especially true when the vehicle is already en route.
In some rare cases, however, researchers enable \emph{en route diversion} (e.g., \citet{ulmer2017value, bosse2023dynamic}).
In some other cases, vehicle movements are simulated within a real-world road network, where turning on the street is not allowed, so diversions from the current route can only take place at the next road crossing (e.g., \citet{ferrucci2014realtime, ferrucci2015general, ferrucci2016proactive}).
Since in these problems, the road crossings can also be modeled as locations, we do not consider this approach as an en route diversion.
In a similar approach, \citet{haferkamp2024design} considered those locations to be deviation points that were located on a traveled shortest path.

\subsubsection{Split delivery}\label{sec:int:dec:split}

\emph{Split delivery} means that a single request can be served by multiple vehicles (or multiple routes of the same vehicle).
Although split delivery is more typical of VRPDSRs (e.g., \citet{schyns2015ant, sarasola2016variable}), it also occurs in some DPDPs.
In the problem formulation of \citet{hao2022introduction} for a DPDP, orders are inherently split into the smallest deliverable units, and can only be shipped separately if their total demand exceeds the uniform vehicle capacity.

\section{A general modeling framework for dynamic vehicle routing I. - Basic concepts}\label{sec:basics}

In this section, we propose the basic concept and terminology of our modeling framework.
First, we provide an overview of the problems under investigation (\Cref{sec:bas:overview}).
Then, we discuss the main elements in detail, which are the locations (\Cref{sec:bas:locations}), the orders (\Cref{sec:bas:orders}), and the vehicles (\Cref{sec:bas:vehicles}).

\subsection{Main overview: modeling scope}\label{sec:bas:overview}

A heterogeneous fleet of vehicles must serve pickup-and-delivery type orders that arrive dynamically in the planning horizon.
The pickup/delivery locations can refer to a designated depot, so our modeling framework is suitable for modeling not only DPDPs, but also SDDPs.
Various VRPDSRs can be modeled, for example, by specifying coincident pickup and delivery locations.
Due to the dynamic nature of the problem, the decision maker has the opportunity to re-plan the vehicle routes at certain decision points.
Decision points may be imposed by arbitrary events (e.g., on order request, on vehicle arrival) or may occur periodically.
Any parameter of the problem (e.g., request of orders, travel times) can be deterministic or stochastic.

A service time window can be associated with both the pickup and the delivery of the orders.
Both cancellation by the customers and rejection by the decision maker can be handled.
In the latter case, the postponement of the decision on acceptance/rejection is also allowed.

Split deliveries are allowed, but in this case, the orders must be split into the smallest deliverable units in advance.
It is the decision maker's responsibility to combine and assign them to vehicles according to the splitting rules.

Vehicles can be capacitated or uncapacitated, and may be subject to loading rules.
Delaying the departure is possible.
The planned routes of the vehicles can be modified during their execution, however, en route diversion is not allowed.
Locations may have limited docking capacity, so the vehicles may have to wait for service.

\paragraph{Simulation vs. Decision making}

Certain aspects of the problem (\Cref{sec:bas:problem}) and the decisions (\Cref{sec:bas:decision}) are not necessarily subject to simulation, but rather to decision making.
For example, earliest service start times must obviously be considered by the simulation (since the vehicles must be kept waiting), but latest service start times are the responsibility of the decision maker.
Therefore, some aspects, such as order due dates or depot deadlines are not discussed in our modeling framework.
However, they can be easily adapted.

\subsection{Locations}\label{sec:bas:locations}

Locations can refer to different places, such as where orders are to be picked up or delivered, where vehicles are initially located, or they can represent intersections in the real road network.
The physical movement of vehicles between locations is not detailed, we just assume that after a vehicle departed for its next location, it will arrive there after a certain amount of time.
This travel time must be given or calculable between any two locations that may appear consecutively in the vehicle's route plan, see later in \Cref{sec:bas:veh:route}.
Travel times can be stochastic.

\subsection{Orders}\label{sec:bas:orders}

Each order~$o_i$ has a \emph{pickup location}~$l_i^p$ and a \emph{delivery location}~$l_i^d$, which can refer to depots.
An order~$o_i$ is requested at its \emph{release time}~$r_i$ (for static orders, if any, $r_i = 0$). 
Orders may be associated with an earliest start time for both pickup and delivery.
If the vehicle arrives early, it must wait until the latest earliest start time.

\subsubsection{Order postponement}\label{sec:bas:ord:post}

In our approach, the decision on an order (i.e., accept/reject) can be postponed until a specific time point.
Assume that a decision is made at time~$t_1$ in which an order is postponed until time~$t_2$.
The postponement means the following.

\paragraph{Case 1 (postponement is expired)}
If no decision point is imposed in time interval $[t_1,t_2]$, the postponement of the order is expired.
Thus, a decision point will be imposed at~$t_2$, which enables the decision maker to reconsider the order.

\paragraph{Case 2 (postponement is interrupted)}
If a decision point is imposed in $[t_1,t_2]$, the postponement of the order will be interrupted at that time.
The decision maker may now accept/reject the order, or postpone it again.

\subsection{Vehicles}\label{sec:bas:vehicles}
We consider a heterogeneous fleet of vehicles, denoted with $\mathcal{V}$.
Each vehicle~$v$ is associated with an \emph{initial location} $l^{\operatorname{init}}_v$.

\subsubsection{Route plans}\label{sec:bas:veh:route}

The movements of the vehicles are controlled by their route plans.
The \emph{route plan} of a vehicle~$v$ is a sequence of visits
\[
  \theta_v = \left( \theta^j_v : j = 1,\ldots,\ell_v \right)\ \text{with}\ 
  \theta^j_v = \left( l^j_v, \mathcal{P}^j_v, \mathcal{D}^j_v; est^j_v \right),
\]
where each \emph{visit}~$\theta^j_v$ is specified by a location ($l^j_v$) to which the vehicle must travel (unless it is currently there), and by (possibly empty) ordered lists ($\mathcal{P}^j_v$ and $\mathcal{D}^j_v$) containing the orders that must be picked up and delivered at the location, respectively.
In addition, an \emph{earliest start time} ($est^j_v$) can be associated with the visit, indicating the earliest time when the vehicle can depart for that location, see later in \Cref{sec:bas:veh:del}.
Route plans will be used later in our decision process to describe the states (\Cref{sec:sdp:state}) and the decisions (\Cref{sec:sdp:actions}).
For an insightful example of route plans we also refer to that section (\Cref{sec:sdp:example}).

\subsubsection{Execution of the route plans}\label{sec:bas:veh:exe}

\begin{figure*}
\centering
\begin{tikzpicture}
\tikzstyle{stateinter}   = [thick,->,latex-latex]
\tikzstyle{timenode}     = [below=0.1,scale=.8,align=center]
\tikzstyle{statenode}    = [circle,fill=black,inner sep=0pt,minimum size= 1ex]
\tikzstyle{movementnode} = [above,scale=.8,align=center]
\tikzstyle{statusnode}   = [above,scale=.8,align=center,font=\it]
\tikzstyle{timedashed}   = [dashed]

\def\xdep{0}
\def\xarr{4}
\def\xser{6}
\def\xfin{9}
\def\xnext{11.25}

\def\levone{0.75}
\def\levtwo{1.75}
\def\lewthr{2.75}

\draw[thick,->,-stealth] (-0.75,0) -- (12,0);

\node[statenode] (sdep) at (\xdep,0) {};
\node[timenode] at (sdep) {departure};

\node[statenode] (sarr) at (\xarr,0) {};
\node[timenode] (sarrlabel) at (sarr) {arrival};

\node[statenode] (sser) at (\xser,0) {};
\node[timenode] (sserlabel) at (sser) {service\\start};

\node[statenode] (sfin) at (\xfin,0) {};
\node[timenode] (sfinlabel) at (sfin) {service\\finish};

\node[statenode] (sdep2) at (\xnext,0) {};
\node[timenode] (sdep2label) at (sdep2) {departure};

\draw[stateinter] (\xdep,\levone) -- (\xarr,\levone);
\node[movementnode] at ($(\xdep,\levone)!0.5!(\xarr,\levone)$) {travel};

\draw[stateinter] (\xarr,\levone) -- (\xser,\levone);
\node[movementnode] at ($(\xarr,\levone)!0.5!(\xser,\levone)$) {pre-service};

\draw[stateinter] (\xser,\levone) -- (\xfin,\levone);
\node[movementnode] at ($(\xser,\levone)!0.5!(\xfin,\levone)$) {service};

\draw[stateinter] (\xfin,\levone) -- (\xnext,\levone);
\node[movementnode] at ($(\xfin,\levone)!0.5!(\xnext,\levone)$) {pre-departure};

\draw[stateinter] (\xarr,\levtwo) -- (\xser,\levtwo);
\node[statusnode] at ($(\xarr,\levtwo)!0.5!(\xser,\levtwo)$) {waiting\\for service};

\draw[stateinter] (\xser,\levtwo) -- (\xfin,\levtwo);
\node[statusnode] at ($(\xser,\levtwo)!0.5!(\xfin,\levtwo)$) {under service};

\draw[stateinter] (\xfin,\levtwo) -- (\xnext,\levtwo);
\node[statusnode] at ($(\xfin,\levtwo)!0.5!(\xnext,\levtwo)$) {idle};

\draw[stateinter] (\xdep,\lewthr) -- (\xarr,\lewthr);
\node[statusnode] at ($(\xdep,\lewthr)!0.5!(\xarr,\lewthr)$) {en route};

\draw[stateinter] (\xarr,\lewthr) -- (\xnext,\lewthr);
\node[statusnode] (atloclabel) at ($(\xarr,\lewthr)!0.5!(\xnext,\lewthr)$) {at location};

\draw[timedashed] (\xdep,0)  -- (\xdep,\lewthr+0.5);
\draw[timedashed] (\xarr,0)  -- (\xarr,\lewthr+0.5);
\draw[timedashed] (\xser,0)  -- (\xser,\levtwo+0.5);
\draw[timedashed] (\xfin,0)  -- (\xfin,\levtwo+0.5);
\draw[timedashed] (\xnext,0) -- (\xnext,\lewthr+0.5);

\draw [decorate,decoration={brace,amplitude=5pt,mirror,raise=5ex}]
  (\xarr,0) -- (\xnext,0) node[midway,yshift=-8ex,scale=.8]{visit};
\end{tikzpicture}
\caption{Vehicle operations between two consecutive departures.}
\label{fig:mod:veh:exe}
\end{figure*}
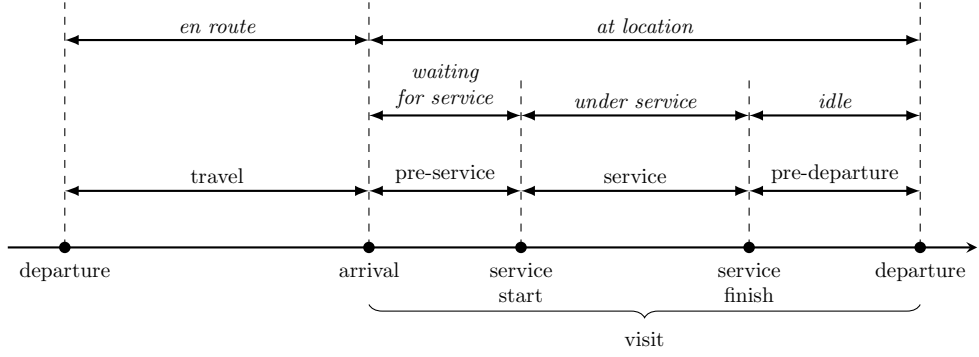

Vehicles -- according to their route plan -- travel from location to location to perform services there, i.e. to pickup and/or deliver orders.
In \Cref{fig:mod:veh:exe}, we depicted the vehicle operations.

\paragraph{Travel}

By \emph{travel}, we mean that the vehicle departs from its \emph{current location} to a specific location, called \emph{destination}.
From \emph{departure} to \emph{arrival}, the vehicle is \emph{en route} (i.e., \emph{on the way}).
While the vehicle is en route, its exact position is not known.
Consequently, the travel cannot be interrupted nor redirected, that is, once the vehicle departed from its current location, it must arrive sooner or later at its destination.

\paragraph{Service}

At locations, vehicles perform services.
The \emph{service} includes the delivery (unloading) and the pickup (loading) of the corresponding orders, if any, but it may also include other operations, for example, parking or docking.
During the service, the vehicle is \emph{under service}.
Note that the service may be void, for example, when empty vehicles return back to a depot, or when the location represents a road crossing.
Similar to travel, the service cannot be interrupted.

\paragraph{Pre-service}

When a vehicle arrives at a location, its service may not start immediately for various reasons.
For example, some orders may have an earliest service start time that has not yet passed, some orders may not be ready upon arrival, or some docking restrictions may delay the service.
The period between the arrival and the subsequent service start is called \emph{pre-service}.
During this period, we say the vehicle is \emph{waiting for service}.

\paragraph{Pre-departure}

When the service is finished, the vehicle may not depart immediately for various reasons.
For example, the vehicle may have completed its route plan, so the vehicle remains at that location until a new route plan is set.
Or the vehicle may have a remaining route, but the start of its execution has been postponed to a later time (see later in \Cref{sec:bas:veh:del}).
The period between the service finish and the subsequent departure is called \emph{pre-departure}.
During this period, we say the vehicle is \emph{idle}.

\subsubsection{Delaying the departure}\label{sec:bas:veh:del}

Now, we describe our concept for delaying the departure of the vehicles.
Assume that vehicle~$v$ is ready to departure at time~$t_1$ to its next location, however, an earliest start time~$t_2$ is associated with its next visit.
Delaying the departure means the following.

\paragraph{Case 1 (departure postponement is expired)}
If no decision point is imposed in time interval $[t_1,t_2]$, then the postponement of the vehicle is expired.

\paragraph{Case 1.1 (decision point on departure postponement expiration)}
If decision points must be imposed on postponement expiration, then a decision point is imposed at~$t_2$, which allows the decision maker, for example, to re-plan the route of the vehicle.

\paragraph{Case 1.2 (no decision point is needed)}
If no decision points need to be imposed on postponement expiration, then the vehicle departs toward its next location to visit. 

\paragraph{Case 2 (departure postponement is interrupted)}
If a decision point is imposed at $[t_1,t_2]$, then the postponement of the vehicle's departure is interrupted at that time.
The decision maker may re-plan the route of the vehicle.

\section{A general modeling framework for dynamic vehicle routing II. - Discrete-event based decision process}
\label{sec:mod}

\begin{figure}
\centering
\begin{tikzpicture}[scale=0.8,transform shape]
\node[fc-startstop] (start) at (0,0) {start};
\node[fc-state,below=0.5cm of start] (initial) {$s\gets$ initial state};
\node[fc-quest,below=1.0cm of initial] (queue) {emtpy queue?};
\node[fc-state,right=1.0cm of queue] (nextevent) {$e\gets$ next event};
\node[fc-state,right=0.5cm and 0.5 of nextevent] (transition) {transition\\$s\gets\phi(s,e)$};
\node[fc-quest,right=0.5cm of transition] (devent) {decision point?};
\node[fc-state,right=1.0cm of devent] (newevents) {adjust queue};
\node[fc-state,below=1.0cm of devent] (dm) {decision making};
\node[fc-startstop,below=1.0cm of queue] (stop) {stop};

\draw[fc-arrow] (start) -- (initial);
\draw[fc-arrow] (initial) -- (queue);
\draw[fc-arrow] (queue) -- node[fc-arrowlabel,midway] {no} (nextevent);
\draw[fc-arrow] (nextevent) -- (transition);
\draw[fc-arrow] (transition) -- (devent);
\draw[fc-arrow] (devent) -- node[fc-arrowlabel,midway] {no} (newevents);
\draw[fc-arrow] (queue) -- node[fc-arrowlabel,midway] {yes} (stop);
\draw[fc-arrow] (newevents.north) -- ++(0,0.5) -| (queue.after north east);
\draw[fc-arrow] (devent) -- node[fc-arrowlabel,midway] {yes} (dm);
\draw[fc-arrow] (dm) -| (newevents);

\scoped[on background layer]
{
\node[fc-fit,fit=(devent)(dm)(newevents)] (eproc) {};
\node[right,above left=0 and 0 of eproc.south east,scale=0.8,font=\it] {event processing};
}
\end{tikzpicture}
\caption{Sketch of the discrete-event based decision process.}
\label{fig:debdp}
\end{figure}
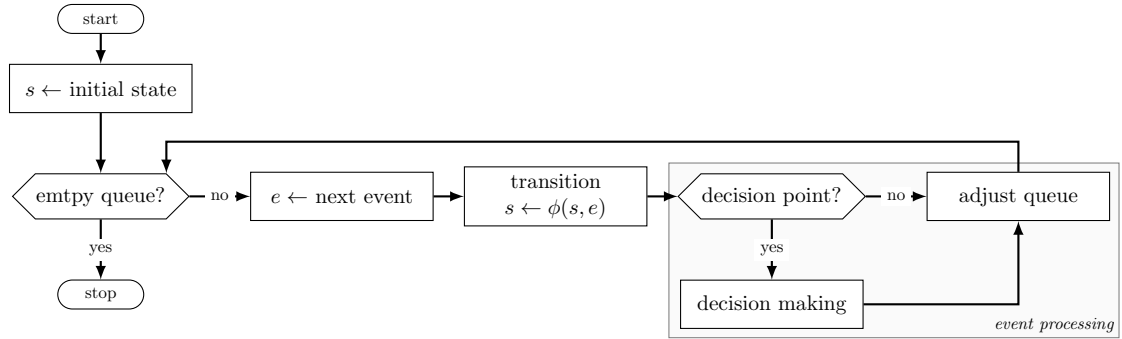

In this section, we propose our modeling framework, which is called \emph{discrete-event based decision process} reflecting on that it is a combination of the discrete-event simulation and the sequential decision process.
The sketch of the process is depicted in \Cref{fig:debdp}.
First, we give a main overview of the framework (\Cref{sec:sdp:main}).
Then, we describe the main elements in detail, which are the states (\Cref{sec:sdp:state}), the events (\Cref{sec:sdp:events}), and the actions (\Cref{sec:sdp:actions}).

\subsection{Main overview}\label{sec:sdp:main}

The status of the system -- including the current position of vehicles and the current status of orders -- is described by \emph{states}.
Various \emph{events} (e.g., an order is requested, a vehicle arrives at a location, etc.) occur in the planning horizon.
These events are stored in an \emph{event queue}, and the decision process jumps from event to event, always to the one associated with the earliest time.
Note that different events can be associated with the same time, and events can be prioritized to establish a processing order between them.
There are two special events, the \emph{decision point} event and the \emph{decision enforcement} event.
When a decision point event occurs, the decision maker is provided with the current state, and then makes a decision that results in an \emph{action}.
This action is set when the corresponding decision enforcement event occurs.

\subsection{States}\label{sec:sdp:state}

A \emph{state} is a tuple
\[
  s = (t_s, \sVehStatus_s, \sOrdStatus_s),
\]
where $t_s$ is the current simulation time, 
$\sVehStatus_s = \{ \sVehStatus_{s,v}: v \in \mathcal{V} \}$ is the status of the vehicles,
and $\sOrdStatus_s$ is the status of the orders, which are discussed in the following.

\subsubsection{Vehicle status}

The status of vehicle~$v$ with respect to state~$s$ is given as a tuple
\[
  \sVehStatus_{s,v} = \left( \sLoad_{s,v}, \theta_{s,v} \right),
\]
where $\sLoad_{s,v}$ is the load, i.e., the list of orders currently carried by the vehicle, and
\[
  \theta_{s,v} = \left( \theta^j_{s,v} : j = 0,\ldots,\ell_{s,v} \right)
\]
is the route plan of the vehicle consisting of a sequence of visits, where
\[
  \theta^0_{s,v} = \left( l^0_{s,v}, \mathcal{P}^0_{s,v}, \mathcal{D}^0_{s,v}; at_{s,v}^0, st_{s,v}^0, ft_{s,v}^0, dt_{s,v}^0 \right)
\]
is the \emph{origin visit}, and
\[
  \theta^j_{s,v} = \left(
  l^j_{s,v},
  \mathcal{P}^j_{s,v},
  \mathcal{D}^j_{s,v};
  est^j_{s,v}
  \right)
  \ \text{for all}\ j = 1,\ldots,\ell_{s,v}.
\]
are the \emph{next visits}.
The origin visit refers to either the current visit of the vehicle, if the vehicle is at a location, or to its previous visit, if the vehicle is en route.
Each visit~$\theta^j_{s,v}$ consists of a location ($l^j_{s,v}$) and two lists of orders to pickup and to deliver ($\mathcal{P}^j_{s,v}$ and $\mathcal{D}^j_{s,v}$), respectively.
The origin visit has four additional elements:
the arrival time ($at_{s,v}^0$),
the service start time ($st_{s,v}^0$),
the service finish time ($ft_{s,v}^0$),
and the departure time ($dt_{s,v}^0$) corresponding to the visit.
The arrival time is always given, but the other times may not be applicable (denoted by $\varnothing$) if the corresponding event has not happened yet.
For example, if the vehicle is currently at a location, then $dt_{s,v}^0 = \varnothing$.
Otherwise, if $dt_{s,v}^0 \neq \varnothing$, the vehicle is currently on the way to its next location~$l^1_{s,v}$.

\subsubsection{Order status}
The status of the orders with respect to state~$s$ is given as a tuple
\[
  \sOrdStatus_s = \left( \mathcal{O}^{\operatorname{open}}_s, \mathcal{O}^{\operatorname{canc}}_s
  \right),
\]
where $\mathcal{O}^{\operatorname{open}}_s$ is the set of \emph{open orders} (i.e., already released, neither canceled nor rejected, and not yet delivered orders),
and $\mathcal{O}^{\operatorname{canc}}_s$ is the set of those orders that are canceled since the last decision point.

\subsubsection{Initial state ($s_0$)}
In the beginning ($t_{s_0} = 0$) vehicles are empty and idle at their initial locations without next visits, i.e., $\sLoad_{s_0,v} = \emptyset$ and $\theta_{s_0,v} = ( ( l^{\operatorname{init}}_v, \emptyset, \emptyset; 0, 0, 0, \varnothing ) )$ for each vehicle~$v$.
No orders are requested yet, that is, $\mathcal{O}^{\operatorname{open}}_{s_0} = \emptyset$ and $\mathcal{O}^{\operatorname{canc}}_{s_0} = \emptyset$.

\subsection{Events}\label{sec:sdp:events}

Each event is associated with a time.
Events are stored in an event queue.
When an event occurs, the state of the system changes (\Cref{sec:sdp:eve:trans}), and then several other events may be inserted to or removed from the event queue (\Cref{sec:sdp:eve:proc}).

Various events can be considered in the model.
In the following (we can call it the default model), we consider the following twelve events: \emph{order request}, \emph{order cancellation}, \emph{order pickup}, \emph{order delivery}, \emph{order postponement expiration}, \emph{vehicle arrival}, \emph{vehicle departure}, \emph{service start}, \emph{service finish}, \emph{departure postponement expiration}, \emph{decision point}, and \emph{decision enforcement}.

The first ten events have a medium priority.
In contrast, decision point events have a high priority, so if multiple events occur at the same time, decision point events are processed last.
In addition, we do not allow multiple decision point events with the same time to be put in the event queue in order to avoid multiple, superfluous decision making.
Decision enforcement events have a low priority, so they are processed before all other events.

\subsubsection{Transition}\label{sec:sdp:eve:trans}

The decision process steps from event to event, and thus the process transitions from state to state.
Formally, \emph{transition} is a function $\phi : \mathcal{S} \times \mathcal{E} \to \mathcal{S}$,
where $\sState$ is the set of all feasible states, and $\mathcal{E}$ is the set of all events.
In fact, only certain events can be considered for a given state (for example, an en route vehicle cannot depart).
For the feasibility of states, see \Cref{sec:apx:feas}.

In the following, we formally define the transition from state~$s_k = (t_k,\sVehStatus_k,\sOrdStatus_k)$ to the subsequent state~$s_{k+1} = (t_{k+1},\sVehStatus_{k+1},\sOrdStatus_{k+1})$ according to event~$e$, i.e., $s_{k+1} = \phi(s_k,e)$.
Since $s_k$ and $s_{k+1}$ differ only in a few parameters, in order to save space, we only indicate the differences between these states.
So first of all, copy the state: $s_{k+1} \gets s_k$.
Regardless of the type of~$e$, $t_{k+1} \gets t$, where $t$ is the time associated with the event.

\paragraph{Order request}

If event~$e$ refers to the request of order~$o_i$, then the order is added to the set of open orders: $\mathcal{O}^{\operatorname{open}}_{k+1} \gets \mathcal{O}^{\operatorname{open}}_k \cup \{o_i\}$.

\paragraph{Order pickup}

If event~$e$ refers to the pickup of order~$o_i$ (i.e., the end of loading) by vehicle~$v$, then the order is added to the carrying order list of the vehicle: $\sLoad_{k+1,v} \gets \sLoad_{k,v} \cup \{o_i\}$.

\paragraph{Order delivery}

If event~$e$ refers to the delivery of order~$o_i$ (i.e., the end of unloading) by vehicle~$v$, then the order is removed from the set of open orders, and from the carrying list of the vehicle: $\mathcal{O}^{\operatorname{open}}_{k+1} \gets \mathcal{O}^{\operatorname{open}}_k \setminus \{o_i\}$ and $\sLoad_{k+1,v} \gets \sLoad_{k,v} \setminus \{o_i\}$.

\paragraph{Order cancellation}

If event~$e$ refers to the cancellation of order~$o_i$, then the order is moved from the set of open orders to the list of canceled orders: $\mathcal{O}^{\operatorname{open}}_{k+1} \gets \mathcal{O}^{\operatorname{open}}_k \setminus \{o_i\}$ and $\mathcal{O}^{\operatorname{canc}}_{k+1} \gets \mathcal{O}^{\operatorname{canc}}_k \cup \{o_i\}$.

\paragraph{Vehicle arrival}

If event~$e$ refers to the arrival of vehicle~$v$, then the origin visit is removed from the route plan: $\theta^0_{k+1,v} \gets (l^1_{k,v}, \mathcal{P}^1_{k,v}, \mathcal{D}^1_{k,v}; t, \varnothing, \varnothing, \varnothing$), $\ell_{k+1,v} \gets \ell_{k,v}-1$, and $\theta^j_{k+1,v} \gets \theta^{j+1}_{k,v}$ for all $j = 1,\ldots,\ell_{k+1,v}$.

\paragraph{Service start}

If event~$e$ refers to the service start of vehicle~$v$, then the service start time of the origin visit is set: $st^0_{k+1} \gets t$.

\paragraph{Service finish}

If event~$e$ refers to the service finish of vehicle~$v$, then the service finish time of the origin visit is set: $ft^0_{k+1} \gets t$.

\paragraph{Vehicle departure}

If event~$e$ refers to the departure of vehicle~$v$, then the departure time of the origin visit is set: $dt^0_{k+1} \gets t$.

\paragraph{Decision enforcement}

If event~$e$ refers to a decision enforcement,
the list of canceled orders is cleared: $\mathcal{O}^{\operatorname{canc}}_{k+1} \gets \emptyset$, and the decision is enforced (see later in \Cref{sec:sdp:act:pds}).

\subsubsection{Event processing}\label{sec:sdp:eve:proc}

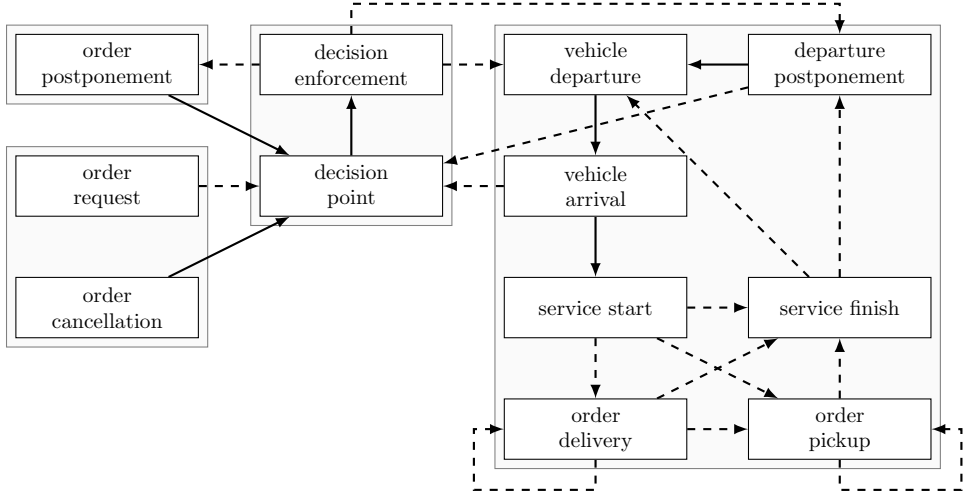
\begin{figure}
\centering
\begin{tikzpicture}[scale=0.8,transform shape]
\tikzstyle{event} = [
  rectangle,
  draw=black,
  fill=white,
  align=center,
  minimum width= 3cm,
  minimum height= 1cm
];

\node[event] (de) {decision\\enforcement};
\node[event,below=1cm of de] (dp) {decision\\point};
\node[event,right=1cm of de] (vd) {vehicle\\departure};
\node[event,below=1cm of vd] (va) {vehicle\\arrival};
\node[event,below=1cm of va] (ss) {service start};
\node[event,right=1cm of ss] (sf) {service finish};
\node[event,right=1cm of vd] (dpe) {departure\\postponement};
\node[event,below=1cm of ss] (od) {order\\delivery};
\node[event,right=1cm of od] (op) {order\\pickup};
\node[event,left=1cm of de] (opp) {order\\postponement};
\node[event,below=1cm of opp] (or) {order\\request};
\node[event,below=1cm of or] (oc) {order\\cancellation};

\scoped[on background layer]{
  \node[fc-fit,fit=(de)(dp)] {};
  \node[fc-fit,fit=(vd)(va)(ss)(od)(op)(sf)(dpe)] {};
  \node[fc-fit,fit=(opp)] {};
  \node[fc-fit,fit=(or)(oc)] {};
  
  \draw[fc-arrow] (vd) -- (va);
  \draw[fc-arrow] (va) -- (ss);
  \draw[fc-arrow,dashed] (ss) -- (sf);
  \draw[fc-arrow,dashed] (ss) -- (od);
  \draw[fc-arrow,dashed] (ss) -- (op);
  \draw[fc-arrow,dashed] (de) -- (vd);
  \draw[fc-arrow,dashed] (od) -- (op);
  \draw[fc-arrow,dashed] (od) -- (sf);
  \draw[fc-arrow,dashed] (op) -- (sf);
  \draw[fc-arrow,dashed] (sf) -- (dpe);
  \draw[fc-arrow,dashed] (sf) -- (vd);
  \draw[fc-arrow] (dpe) -- (vd);
  \draw[fc-arrow,dashed] (od.south) -- ++(0,-0.5) -- ++(-2.0,0.0) |- (od.west);
  \draw[fc-arrow,dashed] (op.south) -- ++(0,-0.5) -- ++(2.0,0.0) |- (op.east);
  \draw[fc-arrow] (dp) -- (de);
  \draw[fc-arrow,dashed] (de) -- (opp);
  \draw[fc-arrow] (opp) -- (dp);
  \draw[fc-arrow,dashed] (or) -- (dp);
  \draw[fc-arrow] (oc) -- (dp);
  \draw[fc-arrow,dashed] (va) -- (dp);
  \draw[fc-arrow,dashed] (dpe) -- (dp);
  \draw[fc-arrow,dashed] (de.north) -- ++(0,0.5) -| (dpe.north);
}
\end{tikzpicture}
\caption{
Events are inductive, meaning that processing one event can cause several new events to be added to or removed from the event queue.
}
\label{fig:events}
\end{figure}

After the transition, the event queue is adjusted, that is, some events may be removed, some new events may be inserted.
In \Cref{fig:events}, we depict which events can induce which other events.
Note that decision points, order request, and order cancellation events can be inserted to the queue from other processes as well.

\paragraph{Decision enforcement}

When a decision enforcement event occurs, the associated action is set (\Cref{sec:sdp:act:pds}).
For each postponed order~$o_i$, if any, an order postponement expired event with time~$\action{pt}_i$ is put into the queue.
For each idle vehicle~$v$, if any, a vehicle departure event with the current time ($\tnow$) or a departure postponement expiration event associated with the earliest start time ($est^1_{s,v}$) is put into the queue, depending on the next visit of the vehicle.

\paragraph{Decision point}

When a decision point event occurs, the decision maker is provided with the current state and returns an action in response.
Then, a decision enforcement event associated with that action and time~$t$ is put into the event queue.
Instantaneous decision making can be modeled with $t=\tnow$ (where $\tnow$ is the current time), while real-time time decision making can be modeled with $t'=\tnow+\delta$, where $\delta$ is the time elapsed during the decision making.
In accordance with \Cref{sec:bas:ord:post,sec:bas:veh:del}, order postponement expiration and departure postponement expiration events, if any, are removed from the queue.

\paragraph{Order requests and cancellations}

When an order request event occurs, a decision point event with time~$\tnow$ may be put into the queue.
On the other hand, if an order cancellation event occurs, it may necessary to insert a decision point event into the queue to prevent the canceled order from being picked up.

\paragraph{Vehicle pre-service}

After the vehicle arrives at a location, a service start event is put into the event queue.
The time associated with the event refers to the time point when the service can be started.
Note that this service start time may depend on the service finish of another vehicles.

\paragraph{Vehicle service}

A vehicle service may consist of several steps.
In the following, we describe the case where orders are first unloaded from the vehicle according to the delivery list, and then orders are loaded to vehicle according to the pickup list.
So, after the service starts, order delivery events, then order pickup events, and finally a service finish event are put into the event queue, one after the other, with the previous one inducing the next.

\paragraph{Vehicle pre-departure}

After the transition triggered by a service finish event, the vehicle can continue to execute its remaining route plan, if any.
(i) If the vehicle has no next visit, there is nothing to do.
(ii) If the vehicle has a next visit, and no earliest start time is associated with it, then a departure event is put into the event queue with the actual simulation time (i.e., the vehicle can depart immediately).
(iii) If an earliest start time is associated with the vehicle's next visit, than a departure postponement expired event is put into the event queue with that time.

\paragraph{Vehicle travel}

After the vehicle departures, a vehicle arrival event -- with the time when the vehicle will arrive -- is put into the event queue.

\subsection{Actions}\label{sec:sdp:actions}

An action is given as a tuple
\[
  x = \left( \action{\sVehStatus}_x, \action{\sOrdStatus}_x \right),
\]
where $\action{\sVehStatus}_x = \{ \action{\sVehStatus}_{x,v} : v \in \mathcal{V} \} $ is set of the updated route plans,
and $\action{\sOrdStatus}_x$ is the decision on orders, which are discussed in the following.

\paragraph{Decision on orders}
The decision on orders is given as a tuple
\[
  \action{\sOrdStatus}_x = \left( 	\action{\mathcal{O}}_{x}^{\operatorname{acc}},
  \action{\mathcal{O}}_{x}^{\operatorname{rej}},
  \action{\mathcal{O}}_{x}^{\operatorname{post}}
  \right),
\]
where $\action{\mathcal{O}}_{x}^{\operatorname{acc}}$,  $\action{\mathcal{O}}_{x}^{\operatorname{rej}}$, and $\action{\mathcal{O}}_{x}^{\operatorname{post}}$ are the set of accepted, rejected, and postponed orders, respectively.
Each postponed order~$o_i \in \action{\mathcal{O}}_{x}^{\operatorname{post}}$ has a time point~$\action{pt}_i$ until the decision on the order is postponed (see \Cref{sec:bas:ord:post}).

\paragraph{Updated route plans}

The updated route plan of a vehicle~$v$ is a sequence of visits
\[
  \action{\sVehStatus}_{x,v} = \left( \action{\theta}_{x,v}^j : j = 0,\ldots,\action{\ell}_{x,v} \right)
\]
with origin visit
\[
  \action{\theta}_{x,v}^0 = \left(
    \action{l}_{x,v}^0,
    \action{\mathcal{P}}_{x,v}^0,
    \action{\mathcal{D}}_{x,v}^0
  \right)
\]
and next visits
\[
  \action{\theta}_{x,v}^j = \left(
    \action{l}_{x,v}^j,
    \action{\mathcal{P}}_{x,v}^j,
    \action{\mathcal{D}}_{x,v}^j;
    \action{est}_{x,v}^j
  \right)\ \text{for all}\ j=1,\ldots,\action{\ell}_{x,v}.
\]
Similarly to the states (\Cref{sec:sdp:state}), each visit~$\action{\theta}_{x,v}^j$ consists of a location ($\action{l}^j_{x,v}$), and pickup and delivery lists ($\action{\mathcal{P}}_{x,v}^j$ and $\action{\mathcal{D}}_{x,v}^j$).
With the exception of the origin visit, each visit can be associated with an earliest start time ($\action{est}_{x,v}^j$).
The origin visit -- more precisely, its pickup and delivery lists -- can be modified until the corresponding service starts.
If no changes have been made to the previous state, the origin visit may not be given (denoted with $\action{\theta}_{x,v}^0 = \varnothing$).

\subsubsection{Transition to post-decision state}\label{sec:sdp:act:pds}

Rejected orders, if any, are removed from the list of open orders: $\mathcal{O}^{\operatorname{open}}_{k+1} \gets \mathcal{O}^{\operatorname{open}}_{k} \setminus \action{\mathcal{O}}_{x}^{\operatorname{rej}}$.
Then, the route plans of the vehicles are updated.
That is, $\theta_{k+1,v}^0 \gets \theta_{k,v}^0$ if $\action{\theta}^0_{x,v} = \varnothing$, otherwise $\theta_{k+1,v}^0 \gets (\action{\theta}^0_{x,v}; at^0_{k,v}, \varnothing, \varnothing, \varnothing)$.
Further, $\ell_{k+1,v} \gets \action{\ell}_{x,v}$ and  $\theta^j_{k+1,v} \gets \action{\theta}^j_{x,v}$ for all $j=1,\ldots,\action{\ell}_{x,v}$.

\subsubsection{Feasibility of actions}\label{sec:sdp:act:feas}

An action~$x$ is \emph{feasible} with respect to state~$s$, if the following constraints are satisfied.
For further feasibility conditions, see \Cref{sec:apx:feas}.

\paragraph{Decision on orders}
Exactly one decision must be made on each order, that is
\[
  \action{\mathcal{O}}_{x}^{\operatorname{acc}} \cup
  \action{\mathcal{O}}_{x}^{\operatorname{rej}} \cup
  \action{\mathcal{O}}_{x}^{\operatorname{post}} =
  \mathcal{O}_{s}^{\operatorname{open}}
\]
such that the sets $\action{\mathcal{O}}_{x}^{\operatorname{acc}}$, $\action{\mathcal{O}}_{x}^{\operatorname{rej}}$, and $\action{\mathcal{O}}_{x}^{\operatorname{post}}$ are pairwise disjunctive.

\paragraph{Origin visit}

The origin visit of a vehicle~$v$ cannot be changed if the service has already started (i.e., the vehicle is either under service or idle or en route).
\[
  st_{s,v}^0 \neq \varnothing \Rightarrow \action{\theta}^0_{x,v} = \varnothing
\]

\paragraph{En route diversion}

If vehicle~$v$ is en route, its destination cannot be changed.
\[
  dt_{s,v}^0 \neq \varnothing \Rightarrow \action{l}^1_{x,v} = l^1_{s,v}
\]

\subsection{Example}\label{sec:sdp:example}

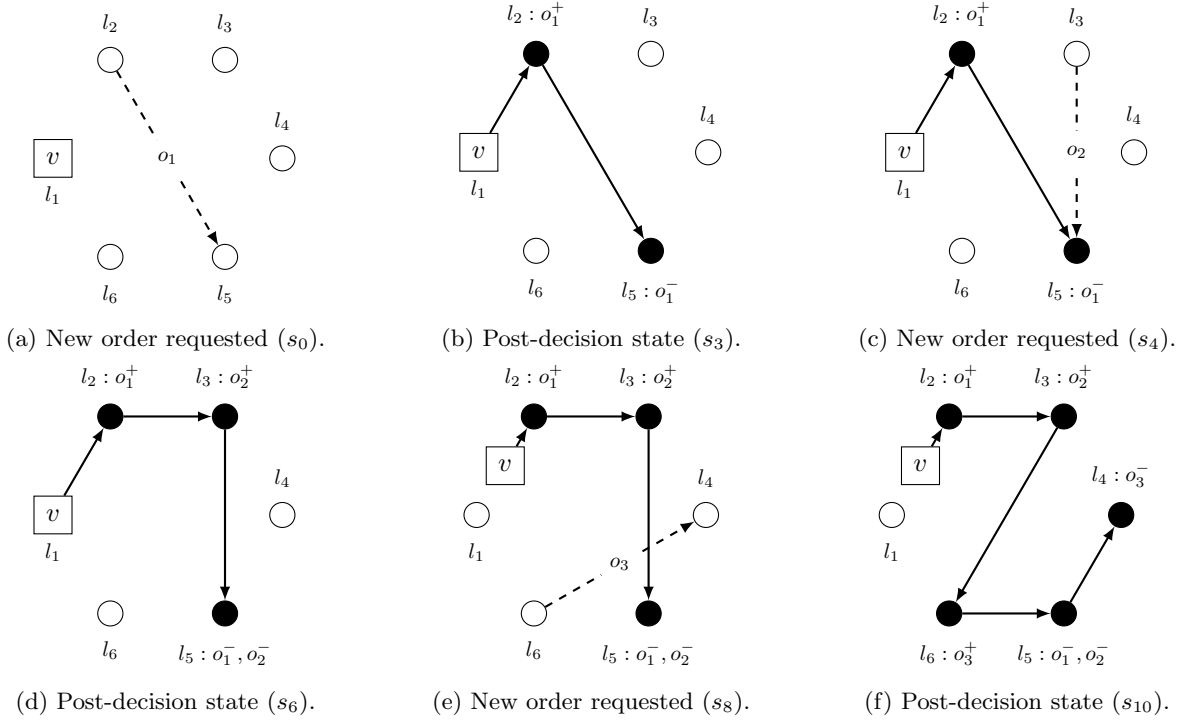
\begin{figure}
\centering
\begin{subfigure}[t]{0.3\textwidth}
\centering
\begin{tikzpicture}
\node[dp-net] (net) {};
\node[dp-loc] (l1) at (net.corner 3) {};
\node[dp-loc] (l2) at (net.corner 2) {};
\node[dp-loc] (l3) at (net.corner 1) {};
\node[dp-loc] (l4) at (net.corner 6) {};
\node[dp-loc] (l5) at (net.corner 5) {};
\node[dp-loc] (l6) at (net.corner 4) {};
\node[dp-veh] (vehicle) at (net.corner 3) {$v$};

\node[dp-loclab,below=0.1cm of l1] {$l_1$};
\node[dp-loclab,above=0.1cm of l2] {$l_2$};
\node[dp-loclab,above=0.1cm of l3] {$l_3$};
\node[dp-loclab,above=0.1cm of l4] {$l_4$};
\node[dp-loclab,below=0.1cm of l5] {$l_5$};
\node[dp-loclab,below=0.1cm of l6] {$l_6$};

\draw[dp-orderarrow] (l2) -- node[dp-arrowlabel] {$o_1$} (l5);
\end{tikzpicture}
\caption{New order requested ($s_0$).}
\end{subfigure}
\hfill%
\begin{subfigure}[t]{0.3\textwidth}
\centering
\begin{tikzpicture}
\node[dp-net]   (net) {};
\node[dp-loc]   (l1) at (net.corner 3) {};
\node[dp-locin] (l2) at (net.corner 2) {};
\node[dp-loc]   (l3) at (net.corner 1) {};
\node[dp-loc]   (l4) at (net.corner 6) {};
\node[dp-locin] (l5) at (net.corner 5) {};
\node[dp-loc]   (l6) at (net.corner 4) {};
\node[dp-veh]   (vehicle) at (net.corner 3) {$v$};

\node[dp-loclab,below=0.1cm of l1] {$l_1$};
\node[dp-loclab,above=0.1cm of l2] {$l_2: o_1^+$};
\node[dp-loclab,above=0.1cm of l3] {$l_3$};
\node[dp-loclab,above=0.1cm of l4] {$l_4$};
\node[dp-loclab,below=0.1cm of l5] {$l_5: o_1^-$};
\node[dp-loclab,below=0.1cm of l6] {$l_6$};

\draw[fc-arrow] (vehicle) -- (l2);
\draw[fc-arrow] (l2) -- (l5);
\end{tikzpicture}
\caption{Post-decision state ($s_3$).}
\end{subfigure}
\hfill%
\begin{subfigure}[t]{0.3\textwidth}
\centering
\begin{tikzpicture}
\node[dp-net]   (net) {};
\node[dp-loc]   (l1) at (net.corner 3) {};
\node[dp-locin] (l2) at (net.corner 2) {};
\node[dp-loc]   (l3) at (net.corner 1) {};
\node[dp-loc]   (l4) at (net.corner 6) {};
\node[dp-locin] (l5) at (net.corner 5) {};
\node[dp-loc]   (l6) at (net.corner 4) {};
\node[dp-veh]   (vehicle) at (net.corner 3) {$v$};

\node[dp-loclab,below=0.1cm of l1] {$l_1$};
\node[dp-loclab,above=0.1cm of l2] {$l_2: o_1^+$};
\node[dp-loclab,above=0.1cm of l3] {$l_3$};
\node[dp-loclab,above=0.1cm of l4] {$l_4$};
\node[dp-loclab,below=0.1cm of l5] {$l_5: o_1^-$};
\node[dp-loclab,below=0.1cm of l6] {$l_6$};

\draw[fc-arrow] (vehicle) -- (l2);
\draw[fc-arrow] (l2) -- (l5);
\draw[dp-orderarrow] (l3) -- node[dp-arrowlabel] {$o_2$} (l5);
\end{tikzpicture}
\caption{New order requested ($s_4$).}
\end{subfigure}
\begin{subfigure}[t]{0.3\textwidth}
\centering
\begin{tikzpicture}
\node[dp-net]   (net) {};
\node[dp-loc]   (l1) at (net.corner 3) {};
\node[dp-locin] (l2) at (net.corner 2) {};
\node[dp-locin] (l3) at (net.corner 1) {};
\node[dp-loc]   (l4) at (net.corner 6) {};
\node[dp-locin] (l5) at (net.corner 5) {};
\node[dp-loc]   (l6) at (net.corner 4) {};
\node[dp-veh]   (vehicle) at (net.corner 3) {$v$};

\node[dp-loclab,below=0.1cm of l1] {$l_1$};
\node[dp-loclab,above=0.1cm of l2] {$l_2: o_1^+$};
\node[dp-loclab,above=0.1cm of l3] {$l_3: o_2^+$};
\node[dp-loclab,above=0.1cm of l4] {$l_4$};
\node[dp-loclab,below=0.1cm of l5] {$l_5: o_1^-, o_2^-$};
\node[dp-loclab,below=0.1cm of l6] {$l_6$};

\draw[dp-arrow] (vehicle) -- (l2);
\draw[dp-arrow] (l2) -- (l3);
\draw[dp-arrow] (l3) -- (l5);
\end{tikzpicture}
\caption{Post-decision state ($s_6$).}
\end{subfigure}
\hfill%
\begin{subfigure}[t]{0.3\textwidth}
\centering
\begin{tikzpicture}
\node[dp-net]   (net) {};
\node[dp-loc]   (l1) at (net.corner 3) {};
\node[dp-locin] (l2) at (net.corner 2) {};
\node[dp-locin] (l3) at (net.corner 1) {};
\node[dp-loc]   (l4) at (net.corner 6) {};
\node[dp-locin] (l5) at (net.corner 5) {};
\node[dp-loc]   (l6) at (net.corner 4) {};
\node[dp-veh]   (vehicle) at (net.side 2) {$v$};

\node[dp-loclab,below=0.1cm of l1] {$l_1$};
\node[dp-loclab,above=0.1cm of l2] {$l_2: o_1^+$};
\node[dp-loclab,above=0.1cm of l3] {$l_3: o_2^+$};
\node[dp-loclab,above=0.1cm of l4] {$l_4$};
\node[dp-loclab,below=0.1cm of l5] {$l_5: o_1^-, o_2^-$};
\node[dp-loclab,below=0.1cm of l6] {$l_6$};

\draw[fc-arrow] (vehicle) -- (l2);
\draw[fc-arrow] (l2) -- (l3);
\draw[fc-arrow] (l3) -- (l5);
\draw[dp-orderarrow] (l6) -- node[dp-arrowlabel] {$o_3$} (l4);
\end{tikzpicture}
\caption{New order requested ($s_8$).}
\end{subfigure}
\hfill%
\begin{subfigure}[t]{0.3\textwidth}
\centering
\begin{tikzpicture}
\node[dp-net]   (net) {};
\node[dp-loc]   (l1) at (net.corner 3) {};
\node[dp-locin] (l2) at (net.corner 2) {};
\node[dp-locin] (l3) at (net.corner 1) {};
\node[dp-locin] (l4) at (net.corner 6) {};
\node[dp-locin] (l5) at (net.corner 5) {};
\node[dp-locin] (l6) at (net.corner 4) {};
\node[dp-veh]   (vehicle) at (net.side 2) {$v$};

\node[dp-loclab,below=0.1cm of l1] {$l_1$};
\node[dp-loclab,above=0.1cm of l2] {$l_2: o_1^+$};
\node[dp-loclab,above=0.1cm of l3] {$l_3: o_2^+$};
\node[dp-loclab,above=0.1cm of l4] {$l_4: o_3^-$};
\node[dp-loclab,below=0.1cm of l5] {$l_5: o_1^-, o_2^-$};
\node[dp-loclab,below=0.1cm of l6] {$l_6: o_3^+$};

\draw[fc-arrow] (vehicle) -- (l2);
\draw[fc-arrow] (l2) -- (l3);
\draw[fc-arrow] (l3) -- (l6);
\draw[fc-arrow] (l6) -- (l5);
\draw[fc-arrow] (l5) -- (l4);
\end{tikzpicture}
\caption{Post-decision state ($s_{10}$).}
\end{subfigure}
\caption{
Selected states from the following scenario:
($s_0$) Vehicle~$v$ is located at location~$l_1$.
($s_1$) Order~$o_1$ is requested.
($s_2$) Decision point is imposed.
($s_3$) Decision maker creates the route plan.
($s_4$) Order~$o_2$ is requested.
($s_5$) Decision point is imposed.
($s_6$) Decision maker updates the route plan.
($s_7$) Vehicle is departed.
($s_8$) Order~$o_2$ is requested.
($s_9$) Decision point is imposed.
($s_{10}$) Decision maker updates the route plan.
}
\label{fig:example}
\end{figure}

In \Cref{fig:example}, we depicted selected states from the following scenario for a dynamic pickup-and-delivery problem.
($s_0$)~Vehicle~$v$ is initially located at location~$l_1$: $t_{s_0} = 0$, $\mathcal{O}_{s_0}^{\operatorname{open}} = \emptyset$, $\theta_{s_0,v} = ((l_1,(),();0,0,0,\varnothing))$.
($s_1$)~Order~$o_1$ from~$l_2$ to~$l_5$ is requested: $\mathcal{O}_{1}^{\operatorname{open}} = \{o_1\}$.
($s_2$)~A decision point is imposed.
The decision maker makes a decision ($x_1$) that the order is accepted, and the initial route plan is created.
However, the departure of the vehicle is delayed until time~10.
That is, $\action{\mathcal{O}}_{x_1}^{\operatorname{acc}} = \{o_1\}$ and
$
\action{\theta}_{x_1,v} = (
  (l_2,(o_1),\emptyset;10),
  (l_5,\emptyset,(o_1);\varnothing)
)
$.
($s_3$)~The decision is enforced:
$
\theta_{s_3,v} = (
  (l_1,\emptyset,\emptyset;0,0,0,\varnothing)$,
  $(l_2,(o_1),();10)$,
  $(l_5,(),(o_1);\varnothing)
)
$.
($s_4$)~Order~$o_2$ from~$l_3$ to~$l_5$ is requested at time~$5$, thus $t_{s_4}=5$, $\mathcal{O}_{s_4}^{\operatorname{open}} = \{o_1,o_2\}$.
($s_5$)~A decision point is imposed.
The decision maker accepts the order and inserts it into the route plan of the vehicle ($x_2$).
That is, $\action{\mathcal{O}}_{x_2}^{\operatorname{acc}} = \{o_1,o_2\}$ and
$
\action{\theta}_{x_2,v} = (
  (l_2,(o_1),\emptyset;10)$,
  $(l_3,(o_2),\emptyset;\varnothing)$,
  $(l_5,\emptyset,(o_1,o_2);\varnothing)
)
$.
($s_6$)~The decision is enforced:
$
\theta_{s_6,v} = (
  (l_1,\emptyset,\emptyset;0,0,0,\varnothing)$,
  $(l_2,(o_1),\emptyset;10)$,
  $(l_3,\emptyset,(o_2);\varnothing)$,
  $(l_5,\emptyset,(o_1,o_2);\varnothing)
)
$.
($s_7$)~The vehicle is departed at time~$10$, that is, $t_{s_7}=10$, and $\theta_{s_7,v}^0 = (l_1,0,0,0,10)$.
($s_8$)~Order~$o_3$ from~$l_4$ to~$l_6$ is requested at time~$12$, that is, $t_{s_8}=12$, and $\mathcal{O}_{s_8}^{\operatorname{open}} = \{o_1,o_2,o_3\}$.
($s_9$)~A decision point is imposed.
The decision maker accepts the order and inserts it into the route plan of the vehicle ($x_3$).
That is, $\action{\mathcal{O}}_{x_3}^{\operatorname{acc}} = \{o_1,o_2,o_3\}$ and
$
\action{\theta}_{x_3,v} = (
  (l_2,(o_1),\emptyset;10)$,
  $(l_3,(o_2),\emptyset;\varnothing)$,
  $(l_6,(o_3),\emptyset;\varnothing)$,
  $(l_5,\emptyset,(o_1,o_2);\varnothing)$,
  $(l_4,\emptyset,(o_3);\varnothing)
)
$.
($s_{10}$)~The decision is enforced:
$
\theta_{s_{10},v} = (
  (l_1,\emptyset,\emptyset;0,0,0,10)$,
  $(l_2,(o_1),\emptyset;10)$,
  $(l_3,(o_2),\emptyset;\varnothing)$,
  $(l_6,(o_3),\emptyset;\varnothing)$,
  $(l_5,\emptyset,(o_1,o_2);\varnothing)$,
  $(l_4,\emptyset,(o_3);\varnothing)
)
$.
($s_{11}$)~The vehicle is arrived at location~$l_2$ at time~$20$: $t_{s_{11}}=20$ and $\theta_{s_{11},v}^0 = (l_2,(o_1),\emptyset;20,\varnothing,\varnothing,\varnothing)$.
($s_{12}$)~After the one-minute parking, the service started: $t_{s_{12}}=21$ and $\theta_{s_{12},v}^0 = (l_2,(o_1),\emptyset;20,21,\varnothing,\varnothing)$.
($s_{13}$)~The loading of order~$o_1$ took two minutes: $t_{s_{13}}=23$, $\sLoad_{s_{13},v} = (o_1)$ and $\theta_{s_{13},v}^0 = (l_2,(o_1),\emptyset;20,21,23,\varnothing)$.
($s_{14}$)~The vehicle departed immediately to its next location: $\theta_{s_{14},v}^0 = (l_2,(o_1),\emptyset;20,21,23,23)$.

\section{An open source simulation tool for dynamic vehicle routing}\label{sec:simulation}

In this section, we briefly present the main components of our simulation framework for dynamic vehicle routing, called \texttt{dvrpsim}.
Our goal is to provide a concise overview of how to use the simulation package.
For an extended, technical description, we refer to \Cref{sec:supp}.
A more detailed tutorial can be found on the webpage of the package: \url{https://sztaki-hu.github.io/dvrpsim/}.

\subsection{A short introduction}

Our simulator is implemented in Python language, however, the implementation of the decision making procedure (also called \emph{external routing algorithm}) is not tied to Python.
For the implementation, we used the SimPy package\footnote{\url{https://simpy.readthedocs.io/en/latest/}}, which is a single-thread process-based discrete-event simulation framework.

\subsubsection{Installation}

The source code is available at \url{https://github.com/sztaki-hu/dvrpsim}.
Assuming Python is already installed, the package can also be installed by typing \texttt{python -m pip install dvrpsim} at the command prompt.

\subsubsection{Modeling (dynamic) vehicle routing problems}

To model a vehicle routing problem, the user needs to build a \texttt{Model}, and to add the necessary \texttt{Location}s, \texttt{Order}s, and \texttt{Vehicle}s that represent the corresponding locations, orders, and vehicles, respectively.
These classes have several callback methods, which can be customized to model their desired behavior.
The routing callback of the \texttt{Model} must be also implemented to connect the external routing algorithm and the simulator.

By starting the simulation (i) each order is requested at its release time; (ii) when a decision point is imposed, the external routing algorithm is called; (iii) once a route plan is set for a vehicle, it begins to execute it.
Unless the user implements otherwise, the simulation ends when all orders have been processed (i.e., delivered, canceled, or rejected).
At the end of the simulation, the history of the vehicles and orders is available, thus various statistics can be generated.

\subsubsection{Locations}

Each \texttt{Location} can optionally be associated with coordinates and a shared resource to model its capacity.
The distances and travel times between the locations can be defined and/or used in the corresponding callbacks of the \texttt{Vehicle}s.

\subsubsection{Orders}

Each \texttt{Order} must be associated with a release time, a pickup location, and a delivery location.
There are also several other optional parameters (such as quantity, pickup/delivery time window, pickup/delivery duration, etc.).

During the simulation, each order is requested at its release time, after which the order is available for insertion into a vehicle route.
Note that orders can also be created on the fly, while the simulation is running.

An \texttt{Order} has several callback methods that are invoked, for example, when the order is requested, rejected, canceled, postponed, picked up, delivered, or when the postponement of the order is expired.
By requesting routing in such a callback, the user can model, for example, decision points on order request/cancellation/postponement.

\subsubsection{Vehicles}

Each \texttt{Vehicle} must be associated with an initial location, and there are several other optional parameters (such as capacity, loading rule, etc.).
In addition, the travel time callback should be defined that returns the travel time for the vehicle between the corresponding locations.

During the simulation, once a route plan is set for a vehicle as a result of decision making, the vehicle begins to execute it.
Recall that the execution procedure of a vehicle consists of four main parts, these are, the pre-departure, the travel, the pre-service, and the service (see \Cref{fig:mod:veh:exe}).
By default, the pre-departure procedure delays the departure of the vehicle when an earliest start time is associated with the next visit.
The travel procedure uses the travel time callback to obtain the arrival time at the next location.
The pre-service procedure takes into account the earliest service start times of the corresponding orders and the capacity of the corresponding location and, if necessary, makes the vehicle wait accordingly.
The service procedure models the unloading and the loading of the corresponding orders.

A \texttt{Vehicle} have several callback methods that are invoked, for example, when the vehicle arrives/departs at/from a location, when the service of the vehicle starts/finishes, or when one of its process is interrupted.
By requesting routing in such a callback, we can model, for example, decision points on vehicle arrival.

\subsubsection{Decision making procedure}

The routing callback of the \texttt{Model} can be used to connect the external routing algorithm and the simulator.
The external routing algorithm can be implemented in arbitrary programming language.
Note that the external routing algorithm does not have to be necessary "external", as the algorithm itself can also be implemented in that callback.

At each decision point, a routing callback is invoked, which includes invoking the external routing algorithm.
The simulator provides the current state in JSON format, allowing file-based interaction with the external routing algorithm, which is especially useful if the latter is not implemented in Python.
The output of the routing algorithm (i.e., the decision) is processed and enforced.
Before enforcing the decision, it is possible to check various problem constraints (e.g., the capacity constraints of the vehicles).
By default, the simulator assumes instantaneous (i.e., zero time) decision making, but real-time decision making can also be modeled.

\subsection{Case studies}

As a proof-of-concept, we implemented several examples using our simulator, which are available together with the source code.
The following three examples deal with three very different problems with very different problem aspects and constraints, demonstrating that the framework is suitable for modeling a wide range of dynamic vehicle routing problems.

\subsubsection{A dynamic pickup-and-delivery problem}

A dynamic pickup-and-delivery problem was introduced in a competition organized by the International Conference on Automated Planning and Scheduling in~2021 (ICAPS 2021), see \citep{hao2022introduction}.

\paragraph{Problem overview}
There is a fleet of homogeneous vehicles that has to serve pickup-and-delivery order requests which occur over a day. Each order is characterized by a quantity, a pickup factory, a delivery factory, a release time, and a due date.
The vehicles can be loaded up to their capacity, while unloading has to follow the last-in-first-out (LIFO) rule.
Those, but only those orders whose quantity exceeds the capacity of the vehicles, can be split and delivered separately.
The travel times and the distances between the factories are given.
Each factory has a given number of docking ports for serving (that is, loading and unloading) the vehicles.
Vehicles are served on a first-come-first-served basis.
If a vehicle arrives at a factory and all ports are occupied, its service cannot begin immediately, but the vehicle has to join the waiting queue.
That is, the vehicle must wait until one of the docking ports becomes free, and no vehicle that arrived earlier is waiting for a port.
The objective is to satisfy all the requests such that a combination of tardiness penalties and traveling distances is minimized.
Decision points occur in every 10 minutes.

\paragraph{Proof-of-concept}
To model this problem, we used the default \texttt{Location} class, where each location is associated with a shared resource to model the  of its docking ports.
We also used the default \texttt{Order} class.
We inherited a custom \texttt{Vehicle} class, where (i)~the travel time callback returns the pre-given travel times; (ii)~the service procedure is extended to model dock approaching of the vehicles.
Capacity and LIFO loading rule are also set for the vehicles.
A pre-defined method is used to impose decision points in every 10 minutes.
The form of states and decisions is also modified, so that the \texttt{Model} can be connected with the already implemented algorithms for the problem.

\subsubsection{A same-day delivery problem}

\citet{voccia2019sameday} introduced a same-day deliver problem for online purchases.
The benchmark instances for their work are publicly available.

\paragraph{Problem overview}
The problem is characterized by a fleet of vehicles operating from a depot and by a set of locations.
Customers request service throughout the day until a fixed cut-off time. Arrivals of requests are described by a known arrival rate and distribution.
Associated with each request is a known service time and a delivery time window at the customer location.
Once requests are made, a vehicle at the depot can be assigned requests and leave the depot immediately.
Alternatively, a vehicle can wait at the depot before being assigned requests.
Once a vehicle leaves the depot, the route for that vehicle is fixed, and the vehicle returns to the depot when it has made all its assigned deliveries.
A request is assigned to a third party when it is no longer feasible for the request to be served by a vehicle at the depot or one of the vehicles en route.
A decison point is imposed as a result of at least one of the following: (i)~a vehicle arrives at the depot; (ii)~a vehicle ends its waiting period; (iii)~a new request arrives and at least one vehicle is waiting at the depot.

\paragraph{Proof-of-concept}
To model this problem, we used the default \texttt{Location} and \texttt{Order} classes.
There is a location for the depot, and there is a separate location for each customer.
Each location is associated with latitude and longitude coordinates.
We inherited a custom \texttt{Vehicle} class, where the travel time callback returns the travel times calculated on Manhattan-distances.
The 'on arrival' callback of the vehicles, and the 'on request' callback of the orders are customized to impose decision points on the appropriate events.
Our demo routing algorithm sets earliest start time for the routes to delay the departure from the depot.

\subsubsection{A restaurant meal delivery problem}

\citet{ulmer2021restaurant} introduced a restaurant meal delivery problem with random ready times.
The benchmark instances for their work are publicly available.

\paragraph{Problem overview}
The problem is characterized by a fleet of vehicles that seeks to fulfill a random set of delivery orders that arrive during the finite order horizon from restaurants located in a service area. Orders occur according to a known stochastic process.
Each realized order is associated with an order time, a delivery location, a pickup restaurant, and a soft deadline.
The time to prepare a customer’s food at each restaurant is random.
Thus, the driver may need to wait for the order’s completion when arriving to a restaurant.
The dispatcher determines which orders are assigned to which vehicles.
Once made, assignments cannot be altered, therefore, assignments can be postponed.
A decision point occurs when a new customer requests service.
A decision point can also be self-imposed, which happens when an order is postponed.

\paragraph{Proof-of-concept}
To model this problem, we used the default \texttt{Location} and \texttt{Order} classes.
There is a separate location for each restaurant, each customer, and each vehicle.
Each location is associated with latitude and longitude coordinates.
Decision points are imposed on order requests.
We inherited a custom \texttt{Vehicle} class, where (i)~the travel time callback returns the travel times calculated on Euclidean-distances; (ii) the pre-service procedure of vehicles are customized to model stochastic ready times.

\section{Conclusion}

In this paper, we focused on developing a simulation tool designed to model a wide range of dynamic vehicle routing problems (DVRPs) to support the dynamic testing of different solution methods.

We began by conducting an extensive literature review to identify the key aspects and common constraints in DVRPs that should be considered in the modeling framework. Based on these findings, we developed a general modeling and simulation framework tailored for simulation purposes.
Finally, we have created an implementation of the framework and made it freely available. As a proof-of-concept, we have implemented several examples
with our framework. These case studies deal with different problems
with very different problem aspects and constraints, demonstrating that the
framework is suitable for modeling a wide range of dynamic vehicle routing problems.

\section*{Acknowledgments}

This research has been supported by the TKP2021-NKTA-01 NRDIO grant on "Research on cooperative production and logistics systems to support a competitive and sustainable economy".
Mark\'o Horv\'ath acknowledges the support of the J\'anos Bolyai Research Scholarship.

\section*{Declaration of Generative AI and AI-assisted technologies in the
writing process}

During the preparation of this work the authors used deepL in order to to check the accuracy of the English text they have created.
After using this tool, the authors reviewed and edited the content as needed and take full responsibility for the content of the publication.


\bibliography{references}

\begin{thebibliography}{101}
\providecommand{\natexlab}[1]{#1}
\providecommand{\url}[1]{\texttt{#1}}
\expandafter\ifx\csname urlstyle\endcsname\relax
  \providecommand{\doi}[1]{doi: #1}\else
  \providecommand{\doi}{doi: \begingroup \urlstyle{rm}\Url}\fi

\bibitem[Ackermann and Rieck(2023)]{ackermann2023novel}
C.~Ackermann and J.~Rieck.
\newblock A novel repositioning approach and analysis for dynamic ride-hailing
  problems.
\newblock \emph{EURO Journal on Transportation and Logistics}, 12:\penalty0
  100109, 2023.

\bibitem[Ackva and Ulmer(2024)]{ackva2024consistent}
C.~Ackva and M.~W. Ulmer.
\newblock Consistent routing for local same-day delivery via micro-hubs.
\newblock \emph{OR Spectrum}, 46\penalty0 (2):\penalty0 375--409, 2024.

\bibitem[Angelelli et~al.(2016)Angelelli, Mansini, and
  Vindigni]{angelelli2016stochastic}
E.~Angelelli, R.~Mansini, and M.~Vindigni.
\newblock The stochastic and dynamic traveling purchaser problem.
\newblock \emph{Transportation Science}, 50\penalty0 (2):\penalty0 642--658,
  2016.
\newblock ISSN 0041-1655.
\newblock \doi{10.1287/trsc.2015.0627}.
\newblock Publisher: {INFORMS}.

\bibitem[Arslan et~al.(2019)Arslan, Agatz, Kroon, and
  Zuidwijk]{arslan2019crowdsourced}
A.~M. Arslan, N.~Agatz, L.~Kroon, and R.~Zuidwijk.
\newblock Crowdsourced delivery—a dynamic pickup and delivery problem with ad
  hoc drivers.
\newblock \emph{Transportation Science}, 53\penalty0 (1):\penalty0 222--235,
  2019.
\newblock ISSN 0041-1655.
\newblock \doi{10.1287/trsc.2017.0803}.

\bibitem[Auad et~al.(2023)Auad, Erera, and Savelsbergh]{auad2023courier}
R.~Auad, A.~Erera, and M.~Savelsbergh.
\newblock Courier satisfaction in rapid delivery systems using dynamic
  operating regions.
\newblock \emph{Omega}, 121:\penalty0 102917, 2023.

\bibitem[Bekta{\c{s}} et~al.(2014)Bekta{\c{s}}, Repoussis, and
  Tarantilis]{bekta2014chapter}
T.~Bekta{\c{s}}, P.~P. Repoussis, and C.~D. Tarantilis.
\newblock Chapter 11: dynamic vehicle routing problems.
\newblock In \emph{Vehicle Routing: Problems, Methods, and Applications, Second
  Edition}, pages 299--347. SIAM, 2014.

\bibitem[Berbeglia et~al.(2010)Berbeglia, Cordeau, and
  Laporte]{berbeglia2010dynamic}
G.~Berbeglia, J.-F. Cordeau, and G.~Laporte.
\newblock Dynamic pickup and delivery problems.
\newblock \emph{European Journal of Operational Research}, 202\penalty0
  (1):\penalty0 8--15, 2010.

\bibitem[Bertsimas et~al.(2019)Bertsimas, Jaillet, and
  Martin]{bertsimas2019online}
D.~Bertsimas, P.~Jaillet, and S.~Martin.
\newblock Online vehicle routing: The edge of optimization in large-scale
  applications.
\newblock \emph{Operations Research}, 67\penalty0 (1):\penalty0 143--162, 2019.
\newblock ISSN 0030-364X.
\newblock \doi{10.1287/opre.2018.1763}.
\newblock Publisher: {INFORMS}.

\bibitem[Bono et~al.(2021)Bono, Dibangoye, Simonin, Matignon, and
  Pereyron]{bono2021solving}
G.~Bono, J.~S. Dibangoye, O.~Simonin, L.~Matignon, and F.~Pereyron.
\newblock Solving multi-agent routing problems using deep attention mechanisms.
\newblock \emph{{IEEE} Transactions on Intelligent Transportation Systems},
  22\penalty0 (12):\penalty0 7804--7813, 2021.
\newblock ISSN 1558-0016.
\newblock \doi{10.1109/TITS.2020.3009289}.
\newblock Conference Name: {IEEE} Transactions on Intelligent Transportation
  Systems.

\bibitem[Bosse et~al.(2023)Bosse, Ulmer, Manni, and Mattfeld]{bosse2023dynamic}
A.~Bosse, M.~W. Ulmer, E.~Manni, and D.~C. Mattfeld.
\newblock Dynamic priority rules for combining on-demand passenger
  transportation and transportation of goods.
\newblock \emph{European Journal of Operational Research}, 309\penalty0
  (1):\penalty0 399--408, 2023.

\bibitem[Braekers et~al.(2016)Braekers, Ramaekers, and
  Van~Nieuwenhuyse]{braekers2016vehicle}
K.~Braekers, K.~Ramaekers, and I.~Van~Nieuwenhuyse.
\newblock The vehicle routing problem: State of the art classification and
  review.
\newblock \emph{Computers \& Industrial Engineering}, 99:\penalty0 300--313,
  2016.

\bibitem[Branke et~al.(2005)Branke, Middendorf, Noeth, and
  Dessouky]{branke2005waiting}
J.~Branke, M.~Middendorf, G.~Noeth, and M.~Dessouky.
\newblock Waiting strategies for dynamic vehicle routing.
\newblock \emph{Transportation science}, 39\penalty0 (3):\penalty0 298--312,
  2005.

\bibitem[Chen et~al.(2022)Chen, Ulmer, and Thomas]{chen2022deep}
X.~Chen, M.~W. Ulmer, and B.~W. Thomas.
\newblock Deep q-learning for same-day delivery with vehicles and drones.
\newblock \emph{European Journal of Operational Research}, 298\penalty0
  (3):\penalty0 939--952, 2022.
\newblock ISSN 0377-2217.
\newblock \doi{10.1016/j.ejor.2021.06.021}.

\bibitem[Chen et~al.(2023)Chen, Wang, Thomas, and Ulmer]{chen2023same}
X.~Chen, T.~Wang, B.~W. Thomas, and M.~W. Ulmer.
\newblock Same-day delivery with fair customer service.
\newblock \emph{European journal of operational research}, 308\penalty0
  (2):\penalty0 738--751, 2023.

\bibitem[Clarke and Wright(1964)]{clarke1964scheduling}
G.~Clarke and J.~W. Wright.
\newblock Scheduling of vehicles from a central depot to a number of delivery
  points.
\newblock \emph{Operations research}, 12\penalty0 (4):\penalty0 568--581, 1964.

\bibitem[C{\^o}t{\'e} et~al.(2023)C{\^o}t{\'e}, de~Queiroz, Gallesi, and
  Iori]{cote2023branch}
J.-F. C{\^o}t{\'e}, T.~A. de~Queiroz, F.~Gallesi, and M.~Iori.
\newblock A branch-and-regret algorithm for the same-day delivery problem.
\newblock \emph{Transportation Research Part E: Logistics and Transportation
  Review}, 177:\penalty0 103226, 2023.

\bibitem[Dantzig and Ramser(1959)]{dantzig1959truck}
G.~B. Dantzig and J.~H. Ramser.
\newblock The truck dispatching problem.
\newblock \emph{Management science}, 6\penalty0 (1):\penalty0 80--91, 1959.

\bibitem[Dayarian and Savelsbergh(2020)]{dayarian2020crowdshipping}
I.~Dayarian and M.~Savelsbergh.
\newblock Crowdshipping and same-day delivery: Employing in-store customers to
  deliver online orders.
\newblock \emph{Production and Operations Management}, 29\penalty0
  (9):\penalty0 2153--2174, 2020.
\newblock ISSN 1937-5956.
\newblock \doi{10.1111/poms.13219}.

\bibitem[Dayarian et~al.(2020)Dayarian, Savelsbergh, and
  Clarke]{dayarian2020sameday}
I.~Dayarian, M.~Savelsbergh, and J.-P. Clarke.
\newblock Same-day delivery with drone resupply.
\newblock \emph{Transportation Science}, 54\penalty0 (1):\penalty0 229--249,
  2020.
\newblock ISSN 0041-1655.
\newblock \doi{10.1287/trsc.2019.0944}.
\newblock Publisher: {INFORMS}.

\bibitem[de~Armas and
  Meli{\'a}n-Batista(2015{\natexlab{a}})]{de2015constrained}
J.~de~Armas and B.~Meli{\'a}n-Batista.
\newblock Constrained dynamic vehicle routing problems with time windows.
\newblock \emph{Soft Computing}, 19:\penalty0 2481--2498, 2015{\natexlab{a}}.

\bibitem[de~Armas and Meli{\'a}n-Batista(2015{\natexlab{b}})]{de2015variable}
J.~de~Armas and B.~Meli{\'a}n-Batista.
\newblock Variable neighborhood search for a dynamic rich vehicle routing
  problem with time windows.
\newblock \emph{Computers \& Industrial Engineering}, 85:\penalty0 120--131,
  2015{\natexlab{b}}.
\newblock ISSN 0360-8352.
\newblock \doi{10.1016/j.cie.2015.03.006}.

\bibitem[Dieter et~al.(2023)Dieter, Stumpe, Ulmer, and
  Schryen]{dieter2023anticipatory}
P.~Dieter, M.~Stumpe, M.~W. Ulmer, and G.~Schryen.
\newblock Anticipatory assignment of passengers to meeting points for
  taxi-ridesharing.
\newblock \emph{Transportation Research Part D: Transport and Environment},
  121:\penalty0 103832, 2023.

\bibitem[Duan et~al.(2020)Duan, Wei, Zhang, and Xia]{duan2020centralized}
L.~Duan, Y.~Wei, J.~Zhang, and Y.~Xia.
\newblock Centralized and decentralized autonomous dispatching strategy for
  dynamic autonomous taxi operation in hybrid request mode.
\newblock \emph{Transportation Research Part C: Emerging Technologies},
  111:\penalty0 397--420, 2020.
\newblock ISSN 0968-090X.
\newblock \doi{10.1016/j.trc.2019.12.020}.

\bibitem[Ehmke and Campbell(2014)]{ehmke2014customer}
J.~F. Ehmke and A.~M. Campbell.
\newblock Customer acceptance mechanisms for home deliveries in metropolitan
  areas.
\newblock \emph{European Journal of Operational Research}, 233\penalty0
  (1):\penalty0 193--207, 2014.
\newblock ISSN 0377-2217.
\newblock \doi{10.1016/j.ejor.2013.08.028}.

\bibitem[Eksioglu et~al.(2009)Eksioglu, Vural, and
  Reisman]{eksioglu2009vehicle}
B.~Eksioglu, A.~V. Vural, and A.~Reisman.
\newblock The vehicle routing problem: A taxonomic review.
\newblock \emph{Computers \& Industrial Engineering}, 57\penalty0 (4):\penalty0
  1472--1483, 2009.

\bibitem[Ferrucci and Bock(2014)]{ferrucci2014realtime}
F.~Ferrucci and S.~Bock.
\newblock Real-time control of express pickup and delivery processes in a
  dynamic environment.
\newblock \emph{Transportation Research Part B: Methodological}, 63:\penalty0
  1--14, 2014.
\newblock ISSN 0191-2615.
\newblock \doi{10.1016/j.trb.2014.02.001}.

\bibitem[Ferrucci and Bock(2015)]{ferrucci2015general}
F.~Ferrucci and S.~Bock.
\newblock A general approach for controlling vehicle en-route diversions in
  dynamic vehicle routing problems.
\newblock \emph{Transportation Research Part B: Methodological}, 77:\penalty0
  76--87, 2015.
\newblock ISSN 0191-2615.
\newblock \doi{10.1016/j.trb.2015.03.003}.

\bibitem[Ferrucci and Bock(2016)]{ferrucci2016proactive}
F.~Ferrucci and S.~Bock.
\newblock Pro-active real-time routing in applications with multiple request
  patterns.
\newblock \emph{European Journal of Operational Research}, 253\penalty0
  (2):\penalty0 356--371, 2016.
\newblock ISSN 0377-2217.
\newblock \doi{10.1016/j.ejor.2016.02.016}.

\bibitem[Ghiani et~al.(2022)Ghiani, Manni, and Manni]{ghiani2022scalable}
G.~Ghiani, A.~Manni, and E.~Manni.
\newblock A scalable anticipatory policy for the dynamic pickup and delivery
  problem.
\newblock \emph{Computers \& Operations Research}, 147:\penalty0 105943, 2022.
\newblock ISSN 0305-0548.
\newblock \doi{10.1016/j.cor.2022.105943}.

\bibitem[Goel et~al.(2019)Goel, Maini, and Bansal]{goel2019vehicle}
R.~Goel, R.~Maini, and S.~Bansal.
\newblock Vehicle routing problem with time windows having stochastic customers
  demands and stochastic service times: Modelling and solution.
\newblock \emph{Journal of Computational Science}, 34:\penalty0 1--10, 2019.

\bibitem[Goodson et~al.(2016)Goodson, Thomas, and
  Ohlmann]{goodson2016restockingbased}
J.~C. Goodson, B.~W. Thomas, and J.~W. Ohlmann.
\newblock Restocking-based rollout policies for the vehicle routing problem
  with stochastic demand and duration limits.
\newblock \emph{Transportation Science}, 50\penalty0 (2):\penalty0 591--607,
  2016.
\newblock ISSN 0041-1655.
\newblock \doi{10.1287/trsc.2015.0591}.
\newblock Publisher: {INFORMS}.

\bibitem[Gy{\"o}rgyi and Kis(2019)]{gyorgyi2019probabilistic}
P.~Gy{\"o}rgyi and T.~Kis.
\newblock A probabilistic approach to pickup and delivery problems with time
  window uncertainty.
\newblock \emph{European Journal of Operational Research}, 274\penalty0
  (3):\penalty0 909--923, 2019.

\bibitem[Haferkamp(2024)]{haferkamp2024design}
J.~Haferkamp.
\newblock Design of multi-optional pickup time offers in ride-sharing systems.
\newblock \emph{EURO Journal on Transportation and Logistics}, page 100134,
  2024.

\bibitem[Haferkamp and Ehmke(2022)]{haferkamp2022effectiveness}
J.~Haferkamp and J.~F. Ehmke.
\newblock Effectiveness of demand and fulfillment control in dynamic fleet
  management of ride-sharing systems.
\newblock \emph{Networks}, 79\penalty0 (3):\penalty0 314--337, 2022.
\newblock ISSN 1097-0037.
\newblock \doi{10.1002/net.22062}.

\bibitem[Haghani and Jung(2005)]{haghani2005dynamic}
A.~Haghani and S.~Jung.
\newblock A dynamic vehicle routing problem with time-dependent travel times.
\newblock \emph{Computers \& operations research}, 32\penalty0 (11):\penalty0
  2959--2986, 2005.

\bibitem[Hao et~al.(2022)Hao, Lu, Li, Tong, Xiang, Yuan, and
  Zhuo]{hao2022introduction}
J.~Hao, J.~Lu, X.~Li, X.~Tong, X.~Xiang, M.~Yuan, and H.~H. Zhuo.
\newblock Introduction to the dynamic pickup and delivery problem benchmark --
  {ICAPS} 2021 competition, 2022.

\bibitem[He et~al.(2019)He, Han, Cheng, Fan, and Dong]{he2019evolutionary}
Z.~He, G.~Han, T.~C.~E. Cheng, B.~Fan, and J.~Dong.
\newblock Evolutionary food quality and location strategies for restaurants in
  competitive online-to-offline food ordering and delivery markets: An
  agent-based approach.
\newblock \emph{International Journal of Production Economics}, 215:\penalty0
  61--72, 2019.
\newblock ISSN 0925-5273.
\newblock \doi{10.1016/j.ijpe.2018.05.008}.

\bibitem[Heitmann et~al.(2023)Heitmann, Soeffker, Ulmer, and
  Mattfeld]{heitmann2023combining}
R.-J.~O. Heitmann, N.~Soeffker, M.~W. Ulmer, and D.~C. Mattfeld.
\newblock Combining value function approximation and multiple scenario approach
  for the effective management of ride-hailing services.
\newblock \emph{EURO Journal on Transportation and Logistics}, 12:\penalty0
  100104, 2023.

\bibitem[Heitmann et~al.(2024)Heitmann, Soeffker, Klawonn, Ulmer, and
  Mattfeld]{heitmann2024accelerating}
R.-J.~O. Heitmann, N.~Soeffker, F.~Klawonn, M.~W. Ulmer, and D.~C. Mattfeld.
\newblock Accelerating value function approximations for dynamic dial-a-ride
  problems via dimensionality reductions.
\newblock \emph{Computers \& Operations Research}, 167:\penalty0 106639, 2024.

\bibitem[Hyland and Mahmassani(2018)]{hyland2018dynamic}
M.~Hyland and H.~S. Mahmassani.
\newblock Dynamic autonomous vehicle fleet operations: Optimization-based
  strategies to assign {AVs} to immediate traveler demand requests.
\newblock \emph{Transportation Research Part C: Emerging Technologies},
  92:\penalty0 278--297, 2018.
\newblock ISSN 0968-090X.
\newblock \doi{10.1016/j.trc.2018.05.003}.

\bibitem[Ichoua et~al.(2006)Ichoua, Gendreau, and Potvin]{ichoua2006exploiting}
S.~Ichoua, M.~Gendreau, and J.-Y. Potvin.
\newblock Exploiting knowledge about future demands for real-time vehicle
  dispatching.
\newblock \emph{Transportation Science}, 40\penalty0 (2):\penalty0 211--225,
  2006.

\bibitem[Jeong and Moon(2024)]{jeong2024dynamic}
J.~Jeong and I.~Moon.
\newblock Dynamic pickup and delivery problem for autonomous delivery robots in
  an airport terminal.
\newblock \emph{Computers \& Industrial Engineering}, 196:\penalty0 110476,
  2024.

\bibitem[Karami et~al.(2020)Karami, Vancroonenburg, and
  Vanden~Berghe]{karami2020periodic}
F.~Karami, W.~Vancroonenburg, and G.~Vanden~Berghe.
\newblock A periodic optimization approach to dynamic pickup and delivery
  problems with time windows.
\newblock \emph{Journal of Scheduling}, 23\penalty0 (6):\penalty0 711--731,
  2020.
\newblock ISSN 1094-6136, 1099-1425.
\newblock \doi{10.1007/s10951-020-00650-x}.

\bibitem[Klapp et~al.(2018{\natexlab{a}})Klapp, Erera, and
  Toriello]{klapp2018dynamic}
M.~A. Klapp, A.~L. Erera, and A.~Toriello.
\newblock The dynamic dispatch waves problem for same-day delivery.
\newblock \emph{European Journal of Operational Research}, 271\penalty0
  (2):\penalty0 519--534, 2018{\natexlab{a}}.
\newblock ISSN 0377-2217.
\newblock \doi{10.1016/j.ejor.2018.05.032}.

\bibitem[Klapp et~al.(2018{\natexlab{b}})Klapp, Erera, and
  Toriello]{klapp2018onedimensional}
M.~A. Klapp, A.~L. Erera, and A.~Toriello.
\newblock The one-dimensional dynamic dispatch waves problem.
\newblock \emph{Transportation Science}, 52\penalty0 (2):\penalty0 402--415,
  2018{\natexlab{b}}.
\newblock ISSN 0041-1655.
\newblock \doi{10.1287/trsc.2016.0682}.
\newblock Publisher: {INFORMS}.

\bibitem[Klapp et~al.(2020)Klapp, Erera, and Toriello]{klapp2020request}
M.~A. Klapp, A.~L. Erera, and A.~Toriello.
\newblock Request acceptance in same-day delivery.
\newblock \emph{Transportation Research Part E: Logistics and Transportation
  Review}, 143:\penalty0 102083, 2020.
\newblock ISSN 1366-5545.
\newblock \doi{10.1016/j.tre.2020.102083}.

\bibitem[Kullman et~al.(2022)Kullman, Cousineau, Goodson, and
  Mendoza]{kullman2022dynamic}
N.~D. Kullman, M.~Cousineau, J.~C. Goodson, and J.~E. Mendoza.
\newblock Dynamic ride-hailing with electric vehicles.
\newblock \emph{Transportation Science}, 56\penalty0 (3):\penalty0 775--794,
  2022.

\bibitem[Lin et~al.(2014)Lin, Choy, Ho, Lam, Pang, and Chin]{lin2014decision}
C.~Lin, K.~L. Choy, G.~T.~S. Ho, H.~Y. Lam, G.~K.~H. Pang, and K.~S. Chin.
\newblock A decision support system for optimizing dynamic courier routing
  operations.
\newblock \emph{Expert Systems with Applications}, 41\penalty0 (15):\penalty0
  6917--6933, 2014.
\newblock ISSN 0957-4174.
\newblock \doi{10.1016/j.eswa.2014.04.036}.

\bibitem[Liu and Luo(2023)]{liu2023demand}
S.~Liu and Z.~Luo.
\newblock On-demand delivery from stores: Dynamic dispatching and routing with
  random demand.
\newblock \emph{Manufacturing \& Service Operations Management}, 25\penalty0
  (2):\penalty0 595--612, 2023.

\bibitem[Liu(2019)]{liu2019optimizationdriven}
Y.~Liu.
\newblock An optimization-driven dynamic vehicle routing algorithm for
  on-demand meal delivery using drones.
\newblock \emph{Computers \& Operations Research}, 111:\penalty0 1--20, 2019.
\newblock ISSN 03050548.
\newblock \doi{10.1016/j.cor.2019.05.024}.

\bibitem[Los et~al.(2020)Los, Schulte, Spaan, and Negenborn]{los2020value}
J.~Los, F.~Schulte, M.~T.~J. Spaan, and R.~R. Negenborn.
\newblock The value of information sharing for platform-based collaborative
  vehicle routing.
\newblock \emph{Transportation Research Part E: Logistics and Transportation
  Review}, 141:\penalty0 102011, 2020.
\newblock ISSN 1366-5545.
\newblock \doi{10.1016/j.tre.2020.102011}.

\bibitem[Ma et~al.(2015)Ma, Zheng, and Wolfson]{ma2015realtime}
S.~Ma, Y.~Zheng, and O.~Wolfson.
\newblock Real-time city-scale taxi ridesharing.
\newblock \emph{{IEEE} Transactions on Knowledge and Data Engineering},
  27\penalty0 (7):\penalty0 1782--1795, 2015.
\newblock ISSN 1558-2191.
\newblock \doi{10.1109/TKDE.2014.2334313}.
\newblock Conference Name: {IEEE} Transactions on Knowledge and Data
  Engineering.

\bibitem[Maciejewski et~al.(2016)Maciejewski, Horni, Nagel, and
  Axhausen]{maciejewski2016dynamic}
M.~Maciejewski, A.~Horni, K.~Nagel, and K.~W. Axhausen.
\newblock Dynamic transport services.
\newblock \emph{The multi-agent transport simulation MATSim}, 23:\penalty0
  145--152, 2016.

\bibitem[Maciejewski et~al.(2017)Maciejewski, Bischoff, H{\"o}rl, and
  Nagel]{maciejewski2017towards}
M.~Maciejewski, J.~Bischoff, S.~H{\"o}rl, and K.~Nagel.
\newblock Towards a testbed for dynamic vehicle routing algorithms.
\newblock In \emph{Highlights of Practical Applications of Cyber-Physical
  Multi-Agent Systems: International Workshops of PAAMS 2017, Porto, Portugal,
  June 21-23, 2017, Proceedings 15}, pages 69--79. Springer, 2017.

\bibitem[Marde{\v{s}}i{\'c} et~al.(2023)Marde{\v{s}}i{\'c}, Erdeli{\'c},
  Cari{\'c}, and Durasevi{\'c}]{mardevsic2023review}
N.~Marde{\v{s}}i{\'c}, T.~Erdeli{\'c}, T.~Cari{\'c}, and M.~Durasevi{\'c}.
\newblock Review of stochastic dynamic vehicle routing in the evolving urban
  logistics environment.
\newblock \emph{Mathematics}, 12\penalty0 (1):\penalty0 28, 2023.

\bibitem[Mitrovi{\'c}-Mini{\'c} and Laporte(2004)]{mitrovic2004waiting}
S.~Mitrovi{\'c}-Mini{\'c} and G.~Laporte.
\newblock Waiting strategies for the dynamic pickup and delivery problem with
  time windows.
\newblock \emph{Transportation Research Part B: Methodological}, 38\penalty0
  (7):\penalty0 635--655, 2004.

\bibitem[Muñoz-Carpintero et~al.(2015)Muñoz-Carpintero, Sáez, Cortés, and
  Núñez]{munoz2015methodology}
D.~Muñoz-Carpintero, D.~Sáez, C.~E. Cortés, and A.~Núñez.
\newblock A methodology based on evolutionary algorithms to solve a dynamic
  pickup and delivery problem under a hybrid predictive control approach.
\newblock \emph{Transportation Science}, 49\penalty0 (2):\penalty0 239--253,
  2015.
\newblock ISSN 0041-1655.
\newblock \doi{10.1287/trsc.2014.0569}.
\newblock Publisher: {INFORMS}.

\bibitem[Ng et~al.(2017)Ng, Lee, Zhang, Wu, and Ho]{ng2017multiple}
K.~K.~H. Ng, C.~K.~M. Lee, S.~Z. Zhang, K.~Wu, and W.~Ho.
\newblock A multiple colonies artificial bee colony algorithm for a capacitated
  vehicle routing problem and re-routing strategies under time-dependent
  traffic congestion.
\newblock \emph{Computers \& Industrial Engineering}, 109:\penalty0 151--168,
  2017.
\newblock ISSN 0360-8352.
\newblock \doi{10.1016/j.cie.2017.05.004}.

\bibitem[Nielsen et~al.(2024)Nielsen, Dahanayaka, Perera, Thibbotuwawa, and
  Kilic]{nielsen2024systematic}
P.~Nielsen, M.~Dahanayaka, H.~N. Perera, A.~Thibbotuwawa, and D.~K. Kilic.
\newblock A systematic review of vehicle routing problems and models in
  multi-echelon distribution networks.
\newblock \emph{Supply Chain Analytics}, page 100072, 2024.

\bibitem[Ninikas and Minis(2014)]{ninikas2014reoptimization}
G.~Ninikas and I.~Minis.
\newblock Reoptimization strategies for a dynamic vehicle routing problem with
  mixed backhauls.
\newblock \emph{Networks}, 64\penalty0 (3):\penalty0 214--231, 2014.
\newblock ISSN 1097-0037.
\newblock \doi{10.1002/net.21567}.

\bibitem[Pillac et~al.(2013)Pillac, Gendreau, Gu{\'e}ret, and
  Medaglia]{pillac2013review}
V.~Pillac, M.~Gendreau, C.~Gu{\'e}ret, and A.~L. Medaglia.
\newblock A review of dynamic vehicle routing problems.
\newblock \emph{European Journal of Operational Research}, 225\penalty0
  (1):\penalty0 1--11, 2013.

\bibitem[Pillac et~al.(2018)Pillac, Guéret, and Medaglia]{pillac2018fast}
V.~Pillac, C.~Guéret, and A.~L. Medaglia.
\newblock A fast reoptimization approach for the dynamic technician routing and
  scheduling problem.
\newblock In L.~Amodeo, E.-G. Talbi, and F.~Yalaoui, editors, \emph{Recent
  Developments in Metaheuristics}, pages 347--367. Springer International
  Publishing, 2018.
\newblock ISBN 978-3-319-58253-5.
\newblock \doi{10.1007/978-3-319-58253-5_20}.

\bibitem[Psaraftis(1980)]{psaraftis1980dynamic}
H.~N. Psaraftis.
\newblock A dynamic programming solution to the single vehicle many-to-many
  immediate request dial-a-ride problem.
\newblock \emph{Transportation Science}, 14\penalty0 (2):\penalty0 130--154,
  1980.

\bibitem[Psaraftis et~al.(2016)Psaraftis, Wen, and
  Kontovas]{psaraftis2016dynamic}
H.~N. Psaraftis, M.~Wen, and C.~A. Kontovas.
\newblock Dynamic vehicle routing problems: Three decades and counting.
\newblock \emph{Networks}, 67\penalty0 (1):\penalty0 3--31, 2016.

\bibitem[Rios et~al.(2021)Rios, Xavier, Miyazawa, Amorim, Curcio, and
  Santos]{rios2021recent}
B.~H.~O. Rios, E.~C. Xavier, F.~K. Miyazawa, P.~Amorim, E.~Curcio, and M.~J.
  Santos.
\newblock Recent dynamic vehicle routing problems: A survey.
\newblock \emph{Computers \& Industrial Engineering}, 160:\penalty0 107604,
  2021.

\bibitem[Ritzinger et~al.(2016)Ritzinger, Puchinger, and
  Hartl]{ritzinger2016survey}
U.~Ritzinger, J.~Puchinger, and R.~F. Hartl.
\newblock A survey on dynamic and stochastic vehicle routing problems.
\newblock \emph{International Journal of Production Research}, 54\penalty0
  (1):\penalty0 215--231, 2016.

\bibitem[Sarasola et~al.(2016)Sarasola, Doerner, Schmid, and
  Alba]{sarasola2016variable}
B.~Sarasola, K.~F. Doerner, V.~Schmid, and E.~Alba.
\newblock Variable neighborhood search for the stochastic and dynamic vehicle
  routing problem.
\newblock \emph{Annals of Operations Research}, 236\penalty0 (2):\penalty0
  425--461, 2016.
\newblock ISSN 1572-9338.
\newblock \doi{10.1007/s10479-015-1949-7}.

\bibitem[Sayarshad and Chow(2015)]{sayarshad2015scalable}
H.~R. Sayarshad and J.~Y.~J. Chow.
\newblock A scalable non-myopic dynamic dial-a-ride and pricing problem.
\newblock \emph{Transportation Research Part B: Methodological}, 81:\penalty0
  539--554, 2015.
\newblock ISSN 0191-2615.
\newblock \doi{10.1016/j.trb.2015.06.008}.

\bibitem[Sayarshad and Oliver~Gao(2018)]{sayarshad2018scalable}
H.~R. Sayarshad and H.~Oliver~Gao.
\newblock A scalable non-myopic dynamic dial-a-ride and pricing problem for
  competitive on-demand mobility systems.
\newblock \emph{Transportation Research Part C: Emerging Technologies},
  91:\penalty0 192--208, 2018.
\newblock ISSN 0968-090X.
\newblock \doi{10.1016/j.trc.2018.04.007}.

\bibitem[Schilde et~al.(2014)Schilde, Doerner, and
  Hartl]{schilde2014integrating}
M.~Schilde, K.~F. Doerner, and R.~F. Hartl.
\newblock Integrating stochastic time-dependent travel speed in solution
  methods for the dynamic dial-a-ride problem.
\newblock \emph{European Journal of Operational Research}, 238\penalty0
  (1):\penalty0 18--30, 2014.
\newblock ISSN 0377-2217.
\newblock \doi{10.1016/j.ejor.2014.03.005}.

\bibitem[Schyns(2015)]{schyns2015ant}
M.~Schyns.
\newblock An ant colony system for responsive dynamic vehicle routing.
\newblock \emph{European Journal of Operational Research}, 245\penalty0
  (3):\penalty0 704--718, 2015.
\newblock ISSN 0377-2217.
\newblock \doi{10.1016/j.ejor.2015.04.009}.

\bibitem[Sluijk et~al.(2023)Sluijk, Florio, Kinable, Dellaert, and
  Van~Woensel]{sluijk2023two}
N.~Sluijk, A.~M. Florio, J.~Kinable, N.~Dellaert, and T.~Van~Woensel.
\newblock Two-echelon vehicle routing problems: A literature review.
\newblock \emph{European Journal of Operational Research}, 304\penalty0
  (3):\penalty0 865--886, 2023.

\bibitem[Soeffker et~al.(2022)Soeffker, Ulmer, and
  Mattfeld]{soeffker2022stochastic}
N.~Soeffker, M.~W. Ulmer, and D.~C. Mattfeld.
\newblock Stochastic dynamic vehicle routing in the light of prescriptive
  analytics: A review.
\newblock \emph{European Journal of Operational Research}, 298\penalty0
  (3):\penalty0 801--820, 2022.

\bibitem[Soeffker et~al.(2024)Soeffker, Ulmer, and
  Mattfeld]{soeffker2024balancing}
N.~Soeffker, M.~W. Ulmer, and D.~C. Mattfeld.
\newblock Balancing resources for dynamic vehicle routing with stochastic
  customer requests.
\newblock \emph{OR Spectrum}, pages 1--43, 2024.

\bibitem[Srour et~al.(2018)Srour, Agatz, and Oppen]{srour2018strategies}
F.~J. Srour, N.~Agatz, and J.~Oppen.
\newblock Strategies for handling temporal uncertainty in pickup and delivery
  problems with time windows.
\newblock \emph{Transportation Science}, 52\penalty0 (1):\penalty0 3--19, 2018.

\bibitem[Steever et~al.(2019)Steever, Karwan, and Murray]{steever2019dynamic}
Z.~Steever, M.~Karwan, and C.~Murray.
\newblock Dynamic courier routing for a food delivery service.
\newblock \emph{Computers \& Operations Research}, 107:\penalty0 173--188,
  2019.
\newblock ISSN 0305-0548.
\newblock \doi{10.1016/j.cor.2019.03.008}.

\bibitem[Tafreshian et~al.(2021)Tafreshian, Abdolmaleki, Masoud, and
  Wang]{tafreshian2021proactive}
A.~Tafreshian, M.~Abdolmaleki, N.~Masoud, and H.~Wang.
\newblock Proactive shuttle dispatching in large-scale dynamic dial-a-ride
  systems.
\newblock \emph{Transportation Research Part B: Methodological}, 150:\penalty0
  227--259, 2021.
\newblock ISSN 0191-2615.
\newblock \doi{10.1016/j.trb.2021.06.002}.

\bibitem[Tirado and Hvattum(2017{\natexlab{a}})]{tirado2017determining}
G.~Tirado and L.~M. Hvattum.
\newblock Determining departure times in dynamic and stochastic maritime
  routing and scheduling problems.
\newblock \emph{Flexible Services and Manufacturing Journal}, 29\penalty0
  (3):\penalty0 553--571, 2017{\natexlab{a}}.
\newblock ISSN 1936-6590.
\newblock \doi{10.1007/s10696-016-9242-x}.

\bibitem[Tirado and Hvattum(2017{\natexlab{b}})]{tirado2017improved}
G.~Tirado and L.~M. Hvattum.
\newblock Improved solutions to dynamic and stochastic maritime pick-up and
  delivery problems using local search.
\newblock \emph{Annals of Operations Research}, 253\penalty0 (2):\penalty0
  825--843, 2017{\natexlab{b}}.
\newblock ISSN 1572-9338.
\newblock \doi{10.1007/s10479-016-2177-5}.

\bibitem[Toth and Vigo(2002)]{toth2002vehicle}
P.~Toth and D.~Vigo.
\newblock \emph{The vehicle routing problem}.
\newblock SIAM, 2002.

\bibitem[Ulmer(2019)]{ulmer2019anticipation}
M.~W. Ulmer.
\newblock Anticipation versus reactive reoptimization for dynamic vehicle
  routing with stochastic requests.
\newblock \emph{Networks}, 73\penalty0 (3):\penalty0 277--291, 2019.
\newblock ISSN 1097-0037.
\newblock \doi{10.1002/net.21861}.

\bibitem[Ulmer(2020)]{ulmer2020dynamic}
M.~W. Ulmer.
\newblock Dynamic pricing and routing for same-day delivery.
\newblock \emph{Transportation Science}, 54\penalty0 (4):\penalty0 1016--1033,
  2020.
\newblock ISSN 0041-1655.
\newblock \doi{10.1287/trsc.2019.0958}.
\newblock Publisher: {INFORMS}.

\bibitem[Ulmer and Streng(2019)]{ulmer2019sameday}
M.~W. Ulmer and S.~Streng.
\newblock Same-day delivery with pickup stations and autonomous vehicles.
\newblock \emph{Computers \& Operations Research}, 108:\penalty0 1--19, 2019.
\newblock ISSN 0305-0548.
\newblock \doi{10.1016/j.cor.2019.03.017}.

\bibitem[Ulmer and Thomas(2018)]{ulmer2018sameday}
M.~W. Ulmer and B.~W. Thomas.
\newblock Same-day delivery with heterogeneous fleets of drones and vehicles.
\newblock \emph{Networks}, 72\penalty0 (4):\penalty0 475--505, 2018.

\bibitem[Ulmer et~al.(2017)Ulmer, Heilig, and Vo{\ss}]{ulmer2017value}
M.~W. Ulmer, L.~Heilig, and S.~Vo{\ss}.
\newblock On the value and challenge of real-time information in dynamic
  dispatching of service vehicles.
\newblock \emph{Business \& Information Systems Engineering}, 59:\penalty0
  161--171, 2017.

\bibitem[Ulmer et~al.(2018)Ulmer, Mattfeld, and K{\"o}ster]{ulmer2018budgeting}
M.~W. Ulmer, D.~C. Mattfeld, and F.~K{\"o}ster.
\newblock Budgeting time for dynamic vehicle routing with stochastic customer
  requests.
\newblock \emph{Transportation Science}, 52\penalty0 (1):\penalty0 20--37,
  2018.

\bibitem[Ulmer et~al.(2019{\natexlab{a}})Ulmer, Goodson, Mattfeld, and
  Hennig]{ulmer2019offlineonline}
M.~W. Ulmer, J.~C. Goodson, D.~C. Mattfeld, and M.~Hennig.
\newblock Offline–online approximate dynamic programming for dynamic vehicle
  routing with stochastic requests.
\newblock \emph{Transportation Science}, 53\penalty0 (1):\penalty0 185--202,
  2019{\natexlab{a}}.
\newblock ISSN 0041-1655.
\newblock \doi{10.1287/trsc.2017.0767}.
\newblock Publisher: {INFORMS}.

\bibitem[Ulmer et~al.(2019{\natexlab{b}})Ulmer, Thomas, and
  Mattfeld]{ulmer2019preemptive}
M.~W. Ulmer, B.~W. Thomas, and D.~C. Mattfeld.
\newblock Preemptive depot returns for dynamic same-day delivery.
\newblock \emph{EURO journal on Transportation and Logistics}, 8\penalty0
  (4):\penalty0 327--361, 2019{\natexlab{b}}.

\bibitem[Ulmer et~al.(2020)Ulmer, Goodson, Mattfeld, and
  Thomas]{ulmer2020modeling}
M.~W. Ulmer, J.~C. Goodson, D.~C. Mattfeld, and B.~W. Thomas.
\newblock On modeling stochastic dynamic vehicle routing problems.
\newblock \emph{EURO Journal on Transportation and Logistics}, 9\penalty0
  (2):\penalty0 100008, 2020.

\bibitem[Ulmer et~al.(2021)Ulmer, Thomas, Campbell, and
  Woyak]{ulmer2021restaurant}
M.~W. Ulmer, B.~W. Thomas, A.~M. Campbell, and N.~Woyak.
\newblock The restaurant meal delivery problem: Dynamic pickup and delivery
  with deadlines and random ready times.
\newblock \emph{Transportation Science}, 55\penalty0 (1):\penalty0 75--100,
  2021.

\bibitem[van Heeswijk et~al.(2019)van Heeswijk, Mes, and
  Schutten]{van2019delivery}
W.~J.~A. van Heeswijk, M.~R.~K. Mes, and J.~M.~J. Schutten.
\newblock The delivery dispatching problem with time windows for urban
  consolidation centers.
\newblock \emph{Transportation Science}, 53\penalty0 (1):\penalty0 203--221,
  2019.
\newblock ISSN 0041-1655.
\newblock \doi{10.1287/trsc.2017.0773}.
\newblock Publisher: {INFORMS}.

\bibitem[Voccia et~al.(2019)Voccia, Campbell, and Thomas]{voccia2019sameday}
S.~A. Voccia, A.~M. Campbell, and B.~W. Thomas.
\newblock The same-day delivery problem for online purchases.
\newblock \emph{Transportation Science}, 53\penalty0 (1):\penalty0 167--184,
  2019.

\bibitem[Vonolfen and Affenzeller(2016)]{vonolfen2016distribution}
S.~Vonolfen and M.~Affenzeller.
\newblock Distribution of waiting time for dynamic pickup and delivery
  problems.
\newblock \emph{Annals of Operations Research}, 236\penalty0 (2):\penalty0
  359--382, 2016.
\newblock ISSN 1572-9338.
\newblock \doi{10.1007/s10479-014-1683-6}.

\bibitem[Wang and Kopfer(2015)]{wang2015rolling}
X.~Wang and H.~Kopfer.
\newblock Rolling horizon planning for a dynamic collaborative routing problem
  with full-truckload pickup and delivery requests.
\newblock \emph{Flexible Services and Manufacturing Journal}, 27\penalty0
  (4):\penalty0 509--533, 2015.
\newblock ISSN 1936-6590.
\newblock \doi{10.1007/s10696-015-9212-8}.

\bibitem[Xiang et~al.(2022)Xiang, Tian, Zhang, Xiao, and
  Jin]{xiang2022pairwise}
X.~Xiang, Y.~Tian, X.~Zhang, J.~Xiao, and Y.~Jin.
\newblock A pairwise proximity learning-based ant colony algorithm for dynamic
  vehicle routing problems.
\newblock \emph{{IEEE} Transactions on Intelligent Transportation Systems},
  23\penalty0 (6):\penalty0 5275--5286, 2022.
\newblock ISSN 1558-0016.
\newblock \doi{10.1109/TITS.2021.3052834}.
\newblock Conference Name: {IEEE} Transactions on Intelligent Transportation
  Systems.

\bibitem[Zhang et~al.(2022)Zhang, Ge, Yang, and Tong]{zhang2022review}
H.~Zhang, H.~Ge, J.~Yang, and Y.~Tong.
\newblock Review of vehicle routing problems: Models, classification and
  solving algorithms.
\newblock \emph{Archives of Computational Methods in Engineering}, 29\penalty0
  (1):\penalty0 195--221, 2022.

\bibitem[Zhang and Van~Woensel(2023)]{zhang2023dynamic}
J.~Zhang and T.~Van~Woensel.
\newblock Dynamic vehicle routing with random requests: A literature review.
\newblock \emph{International Journal of Production Economics}, 256:\penalty0
  108751, 2023.

\bibitem[Zhang et~al.(2023)Zhang, Luo, Florio, and
  Van~Woensel]{zhang2023solving}
J.~Zhang, K.~Luo, A.~M. Florio, and T.~Van~Woensel.
\newblock Solving large-scale dynamic vehicle routing problems with stochastic
  requests.
\newblock \emph{European Journal of Operational Research}, 306\penalty0
  (2):\penalty0 596--614, 2023.
\newblock ISSN 0377-2217.
\newblock \doi{10.1016/j.ejor.2022.07.015}.

\bibitem[Zhang et~al.(2018)Zhang, Ohlmann, and Thomas]{zhang2018dynamic}
S.~Zhang, J.~W. Ohlmann, and B.~W. Thomas.
\newblock Dynamic orienteering on a network of queues.
\newblock \emph{Transportation Science}, 2018.
\newblock \doi{10.1287/trsc.2017.0761}.
\newblock Publisher: {INFORMS}.

\bibitem[Zolfagharinia and Haughton(2014)]{zolfagharinia2014benefit}
H.~Zolfagharinia and M.~Haughton.
\newblock The benefit of advance load information for truckload carriers.
\newblock \emph{Transportation Research Part E: Logistics and Transportation
  Review}, 70:\penalty0 34--54, 2014.
\newblock ISSN 1366-5545.
\newblock \doi{10.1016/j.tre.2014.06.012}.

\bibitem[Zolfagharinia and Haughton(2016)]{zolfagharinia2016effective}
H.~Zolfagharinia and M.~Haughton.
\newblock Effective truckload dispatch decision methods with incomplete advance
  load information.
\newblock \emph{European Journal of Operational Research}, 252\penalty0
  (1):\penalty0 103--121, 2016.
\newblock ISSN 0377-2217.
\newblock \doi{10.1016/j.ejor.2016.01.006}.

\end{thebibliography}

\appendix

\section{Tables for literature review}

In \Cref{tab:lit:vrpdsr,tab:lit:dpdp,tab:lit:sddp} we compiled the reviewed papers.
Abbreviations stand for the following.
Vehicles (VEH): Single (1), Homogeneous fleet (Ho), Heterogeneous fleet (He).
Capacitated vehicles (CAP).
Order time-windows (TW): Soft (S), Hard (H).
Order cancellation (CAN).
Decision points (DPs): Periodic (P), Order request (OR), Vehicle arrival (VA), Self-imposed (SI), New information (NI), Order modification (OM).
Delaying the departure (DEL).
Order rejection (REJ).
Decision postponement (PP).
En route diversion (ERD).

\begin{table}
\scriptsize
\caption{Problem and decision making aspects for DPDPs.}
\label{tab:lit:dpdp}
\begin{tabular}{lccccccccc}
\toprule
paper & VEH & CAP & TW & CAN & DPs & DEL & REJ & PP & ERD \\
\midrule
\citet{ferrucci2014realtime} & He & yes & S & - & P & - & yes & - & (yes) \\
\citet{schilde2014integrating} & Ho & yes & S & - & P & yes & - & - & - \\
\citet{zolfagharinia2014benefit} & He & yes & H & - & P & - & yes & - & (yes) \\
\citet{ma2015realtime} & He & yes & H & - & OR & - & yes & - & - \\
\citet{munoz2015methodology} & He & yes & - & - & OR & - & - & - & - \\
\citet{sayarshad2015scalable} & He & yes & - & yes & OR & - & - & - & - \\
\citet{wang2015rolling} & He & yes & H & - & OR/P & - & yes & - & - \\
\citet{vonolfen2016distribution} & Ho & yes & H & - & OR & yes & - & - & - \\
\citet{zolfagharinia2016effective} & He & yes & H & - & P & yes & yes & - & - \\
\citet{tirado2017determining} & He & yes & H & - & OR, VA & yes & yes & yes & - \\
\citet{tirado2017improved} & He & yes & H & - & OR, VA & - & yes & yes & - \\
\citet{hyland2018dynamic} & Ho & yes & - & - & P & - & - & - & - \\
\citet{sayarshad2018scalable} & He & yes & - & - & OR & - & - & - & - \\
\citet{srour2018strategies} & Ho & yes & H & - & OM & yes & yes & - & - \\
\citet{arslan2019crowdsourced} & He & yes & H & - & NI & - & - & - & - \\
\citet{bertsimas2019online} & Ho & yes & H & - & P & - & yes & yes & - \\
\citet{gyorgyi2019probabilistic} & Ho & yes & H & - & OM & yes & yes & - & - \\
\citet{he2019evolutionary} & Ho & yes & S & - & OR & - & - & - & - \\
\citet{liu2019optimizationdriven} & He & yes & - & - & P & - & - & yes & - \\
\citet{steever2019dynamic} & He & yes & S & - & OR & - & - & - & - \\
\citet{duan2020centralized} & Ho & yes & H & - & P & - & yes & - & - \\
\citet{karami2020periodic} & Ho & - & S & - & P & - & - & - & - \\
\citet{los2020value} & He & yes & H & yes & OR & yes & yes & - & - \\
\citet{tafreshian2021proactive} & Ho & yes & H & - & P & yes & yes & - & - \\
\citet{ulmer2021restaurant} & Ho & - & S & - & OR, SI & - & - & yes & - \\
\citet{ghiani2022scalable} & Ho & - & - & - & OR & - & - & - & - \\
\citet{haferkamp2022effectiveness} & Ho & - & H & - & OR & - & yes & - & - \\
\citet{kullman2022dynamic} & Ho & - & H & - & OR, VA & - & yes & - & - \\
\citet{hao2022introduction} & Ho & yes & S & - & P & - & - & - & - \\
\citet{ackermann2023novel} & Ho & yes & - & - & OR, VA, SI & - & yes & yes & - \\
\citet{auad2023courier} & Ho & yes & S & - & P & - & - & yes & - \\
\citet{bosse2023dynamic} & Ho & yes & - & - & OR & - & yes & - & yes \\
\citet{dieter2023anticipatory} & Ho & yes & - & - & OR & - & - & - & - \\
\citet{heitmann2023combining} & Ho & yes & H & - & OR & - & yes & - & - \\
\citet{ackva2024consistent} & Ho & yes & H & - & OR & yes & yes & - & - \\
\citet{jeong2024dynamic} & Ho & yes & - & - & OR & - & - & - & - \\
\citet{heitmann2024accelerating} & Ho & yes & H & - & OR & - & yes & - & - \\
\citet{haferkamp2024design} & Ho & yes & H & (yes) & OR & - & - & - & (yes) \\
\bottomrule
\end{tabular}
\end{table}

\begin{table}
\scriptsize
\caption{Problem and decision making aspects for SDDPs.}
\label{tab:lit:sddp}
\begin{tabular}{lccccccccc}
\toprule
paper & VEH & CAP & TW & CAN & DPs & DEL & REJ & PP & ERD \\
\midrule
\citet{ehmke2014customer} & Ho & - & H & - & OR & - & yes & - & - \\
\citet{klapp2018dynamic} & 1 & - & - & - & P & - & yes & yes & - \\
\citet{klapp2018onedimensional} & 1 & - & - & - & P & yes & yes & - & - \\
\citet{ulmer2018sameday} & He & yes & H & - & OR & - & yes & - & - \\
\citet{ulmer2019sameday} & Ho & yes & - & - & P & - & - & yes & - \\
\citet{ulmer2019preemptive} & 1 & - & - & - & VA & - & yes & - & - \\
\citet{van2019delivery} & Ho & yes & H & - & P & - & - & yes & - \\
\citet{voccia2019sameday} & Ho & - & H & - & OR, VA, SI & yes & yes & yes & - \\
\citet{dayarian2020crowdshipping} & He & yes & S & - & P, VA & yes & - & yes & - \\
\citet{dayarian2020sameday} & He & yes & H & - & VA & yes & yes & - & - \\
\citet{klapp2020request} & Ho & - & - & - & P, OR & yes & yes & - & - \\
\citet{ulmer2020dynamic} & Ho & - & H & - & OR & - & yes & - & - \\
\citet{chen2022deep} & He & yes & H & - & OR & - & yes & - & - \\
\citet{chen2023same} & He & - & H & - & OR & - & yes & - & - \\
\citet{cote2023branch} & Ho & - & H & - & OR, VA, SI & yes & yes & - & - \\
\citet{liu2023demand} & Ho & yes & H & - & P & - & - & (yes) & - \\
\bottomrule
\end{tabular}
\end{table}

\begin{table}
\scriptsize
\caption{Problem and decision making aspects for VRPDSRs.}
\label{tab:lit:vrpdsr}
\begin{tabular}{lccccccccc}
\toprule
paper & VEH & CAP & TW & CAN & DPs & DEL & REJ & PP & ERD \\
\midrule
\citet{lin2014decision} & Ho & yes & H & yes & OR & - & - & - & - \\
\citet{ninikas2014reoptimization} & Ho & yes & H & - & OR & - & - & - & - \\
\citet{ferrucci2015general} & Ho & - & S & - & P & - & - & - & (yes) \\
\citet{de2015variable} & He & yes & S & - & OR & - & yes & - & - \\
\citet{schyns2015ant} & He & yes & H & yes & NI & - & - & - & - \\
\citet{ferrucci2016proactive} & Ho & - & S & - & P & yes & - & - & (yes) \\
\citet{sarasola2016variable} & Ho & yes & - & - & P & - & - & - & - \\
\citet{angelelli2016stochastic} & 1 & - & - & - & VA & - & - & - & - \\
\citet{goodson2016restockingbased} & Ho & yes & - & - & VA & - & - & - & - \\
\citet{ng2017multiple} & Ho & yes & - & - & VA & - & - & - & - \\
\citet{ulmer2017value} & 1 & - & - & - & OR, VA, SI & - & yes & yes & yes \\
\citet{pillac2018fast} & He & - & H & - & OR,VA & yes & yes & - & - \\
\citet{ulmer2018budgeting} & 1 & - & - & - & VA & yes & yes & - & - \\
\citet{zhang2018dynamic} & 1 & - & H & - & VA, SI & - & yes & yes & - \\
\citet{ulmer2019anticipation} & 1 & - & - & - & VA & - & yes & - & - \\
\citet{ulmer2019offlineonline} & 1 & - & - & - & VA & (yes) & yes & - & - \\
\citet{bono2021solving} & Ho & yes & S & - & VA & - & - & - & - \\
\citet{xiang2022pairwise} & Ho & yes & - & - & P & - & - & - & - \\
\citet{zhang2023solving} & Ho & - & - & - & OR & - & yes & - & - \\
\citet{soeffker2024balancing} & Ho & - & - & - & OR & - & yes & - & - \\
\bottomrule
\end{tabular}
\end{table}

\section{Feasibility of states and actions}\label{sec:apx:feas}

A state~$s$ is \emph{feasible} if the following constraints are satisfied.
An action~$x$ is \emph{feasible with respect to} the feasible state~$s$ if, in addition to the constraints described in \Cref{sec:sdp:act:feas}, the post-decision state $\phi(s,x)$ is feasible.

\subsection{General constraints}

Regardless of the problem, the following constraints must always be taken into account.

\paragraph{Assigned orders}

Only open orders can be assigned to vehicles.
\[
  \bigcup_{v\in\mathcal{V}} \left(
    \sLoad_{s,v} \cup \bigcup_{j=0}^{\ell_{s,v}} \left(
      \mathcal{P}^j_{s,v} \cup \mathcal{D}^j_{s,v} 
    \right)
  \right) \subseteq \mathcal{O}^{\operatorname{open}}_s
\]

\paragraph{Pickup and delivery locations}

Orders can only be picked up at their pickup location ($l^p_{\cdot}$), and can only be delivered at their delivery location ($l^d_{\cdot}$).
\[
  o_i \in \mathcal{P}^j_{s,v} \Rightarrow l^j_{s,v} = l^p_i
\]
\[
  o_i \in \mathcal{D}^j_{s,v} \Rightarrow l^j_{s,v} = l^d_i
\]

\paragraph{Pickup and delivery with the same vehicle}

Orders must be delivered by the same vehicle that picked them up.
\[
o_i \in \mathcal{D}^j_{s,v} \Rightarrow o_i \in \sLoad_{s,v} \cup \bigcup_{k=0}^{j-1} \mathcal{P}^k_{s,v}
\]

\paragraph{Pickup and deliver only once}

Orders can only be picked up and delivered once.
That is, the sets $\sLoad_{s,v}$ and $\mathcal{P}^j_{s,v}$ ($j=0,\ldots,\ell_{s,v}$) must be pairwise disjunctive for each vehicle~$v$.
Similarly, for each vehicle~$v$, the sets $\mathcal{D}^j_{s,v}$ ($j=0,\ldots,\ell_{s,v}$) must be  pairwise disjunctive.

\subsection{Problem specific constraints}

There may be several other constraints for a particular problem at hand (e.g., capacity constraints, loading rules).

\paragraph{Capacity constraints}

If vehicle~$v$ is capacitated, then the total quantity of the loaded orders cannot exceed its capacity~$Q_v$
That is,
\[
\sum_{o_i \in \sLoad_{s,v}} q_i + \sum_{j=0}^{j'} \left( \sum_{o_i \in \mathcal{D}^j_{s,v}} q_i - \sum_{o_i \in \mathcal{P}^j_{s,v}} q_i \right) \leq Q_v \quad \text{for all}\ j' = 0,\ldots,\ell_{s,v},
\]
where it is assumed that unloading takes place first and then loading takes place afterwards.

\section{Supplementary material: A general simulation framework for dynamic vehicle routing}\label{sec:supp}

In this section, we take a closer look at the main components of our simulation framework.

\subsection{Process-based discrete-event simulation}\label{sec:sim:pbdes}
For the implementation of the simulator, we used the SimPy package\footnote{\url{https://simpy.readthedocs.io/en/latest/}}, which is a single-thread process-based discrete-event simulation framework.
This framework works with \texttt{Process}es, which can interact with each other via \texttt{Event}s.
A \texttt{Process} is basically a generator function.
When a \texttt{Process} \emph{yields} an \texttt{Event}, the \texttt{Process} is suspended, and the \texttt{Event} is scheduled, i.e., it is added to the event queue.
The \texttt{Process} is resumed, when this \texttt{Event} is \emph{processed}, i.e., it occurs.

In the following, we will use several types of \texttt{Event}s.
The \texttt{Event} $\operatorname{event}(t)$  will be processed at time~$t$.
The timeout \texttt{Event} $\operatorname{timeout}(\delta)$ will be processed $\delta$ times later after it is scheduled.
The latter is equivalent with $\operatorname{event}(\tnow + \delta)$, where $\tnow$ refers to the current simulation time.
The \texttt{Event} $\operatorname{allof}(e_1,e_2,\ldots)$ will be processed when all of the \texttt{Event}s $e_1,e_2,\ldots$ are processed.

Suspended \texttt{Process}es can be \emph{interrupted}.
When a \texttt{Process} is interrupted, the yielded \texttt{Event} is removed from the event queue, and the underlying function is terminated.

\texttt{Process}es can be also interact with each other via \emph{shared resources}.
For example, a \texttt{Resource} is conceptually a semaphore with a given capacity.
A \emph{request} on a \texttt{Resource} yields an \texttt{Event} which will be processed, when a capacity is freed.

\subsection{Routing procedure}\label{sec:sim:routing}

\begin{algorithm}
\caption{Routing procedure modeling instantaneous decision making}\label{alg:sim:routing}
\begin{algorithmic}[1]
\Input $\varphi$: global boolean variable
\State $\varphi \gets \texttt{False}$
\Yield $\operatorname{timeout}(0)$ with high priority \Comment{decision point}
\If{$\varphi$ is \texttt{True}}
\State\Return
\EndIf
\State $\varphi \gets \texttt{True}$
\State invoke callback on routing start
\State invoke external routing algorithm
\Yield $\operatorname{timeout}(0)$ with low priority \Comment{decision enforcement}
\State invoke callback on routing finish
\State enforce decision
\end{algorithmic}
\end{algorithm}

As described in the previous sections, certain events can impose a decision point, and then the decision maker makes a decision, which can, for example, change the route plans.
In this environment, this procedure is implemented in the \emph{routing procedure}, see \Cref{alg:sim:routing}, and we say that some events request routing.
Once a routing is requested, the routing procedure is added to the environment in a separate \texttt{Process}.
The routing procedure consists of the following steps.

\paragraph{Step 1 (Decision point)}
A timeout \texttt{Event} -- which corresponds to a decision point event (see \Cref{sec:sdp:events}) -- is yielded with zero delay and high priority.
Because of this, those \texttt{Event}s that also occur at that time, will be processed before this \texttt{Event}.
Second, we can detect if multiple routing requests arrive at the same time (i.e., multiple routing \texttt{Process}es are added to the environment), thus we can achieve that only the first can proceed to Step~2, see flag~$\varphi$ in \Cref{alg:sim:routing}.

\paragraph{Step 2 (Callback 'on routing start')}
A callback on routing start is invoked.
By default, this procedure interrupts the postponement \texttt{Process}es of the orders and the vehicles, see later in \Cref{sec:sim:ord:opp,sec:sim:veh:predep}.

\paragraph{Step 3 (Routing)}
The external routing algorithm is invoked.
The external routing algorithm is provided with the current state of the model, then returns a decision on the vehicles and the orders.
Finally, a timeout \texttt{Event} $\operatorname{timeout}(\delta)$ -- which corresponds to a decision enforcement event described in \Cref{sec:sdp:events} -- is yielded to indicate the end of the routing.
By default, $\delta = 0$ regardless of the real run-time of the external routing algorithm, that is, the simulator models instantaneous decisions.
If one wants to model real-time decisions, the delay must correspond to the run-time of the external routing algorithm.

\paragraph{Step 4 (Callback 'on routing finish')}
After the routing is finished, a callback is invoked, which can be used, for example, to check the feasibility of the obtained action (see \Cref{sec:apx:feas}).

\paragraph{Step 5 (Decision enforcement)}
The resulted decision is enforced, as described in \Cref{sec:sdp:act:pds}.
If any state-feasibility constraints (\Cref{sec:sdp:act:feas}) is violated, the simulator terminates with error.
After the status of the orders are updated, the postponement procedure of postponed orders, if any, are started in separate \texttt{Process}es, see later in \Cref{sec:sim:ord:opp}.
After the routes are updated, the execution procedure of idle vehicles, if any, are started in separate \texttt{Process}es, see later in \Cref{sec:sim:veh:predep}.

\subsection{Custom processes}

Arbitrary processes can be added to the simulation environment.
For example, the procedure depicted in \Cref{alg:sim:periodic} can be added to impose decision points at given time steps.

\begin{algorithm}
\caption{Procedure for periodic decision points}\label{alg:sim:periodic}
\begin{algorithmic}[1]
\Input $\Delta$: epoch length
\While{some not rejected orders are not delivered}
\State request routing
\Yield $\operatorname{timeout}(\Delta)$
\EndWhile
\end{algorithmic}
\end{algorithm}

\subsection{Orders}

In the simulator, each order is modeled with an instance of the \texttt{Order} class.

\subsubsection{Requesting orders}

There are multiple ways in the simulator to request orders.
For an example, in the procedure depicted in \Cref{alg:ord:or}, the orders are sorted in ascending order based on their release time.
Then, for each order an \texttt{Event} is yielded that corresponds to an order request event described in \Cref{sec:sdp:events}.

\begin{algorithm}
\caption{An order requesting procedure}\label{alg:ord:or}
\begin{algorithmic}[1]
\State $\mathcal{O} \gets $ orders sorted by their release time ($r_{\cdot}$)
\For{order $o_i \in \mathcal{O}$}
  \State $\delta \gets r_i - \tnow$
  \Yield $\operatorname{timeout}(\delta)$ \Comment{order request}
\EndFor
\end{algorithmic}
\end{algorithm}

\subsubsection{Order postponement}\label{sec:sim:ord:opp}

As described in \Cref{sec:bas:ord:post}, our framework allows order postponement.
Recall that when a decision is enforced, a postponement \texttt{Process} is started separately for each postponed order, if any (see Step~5 in \Cref{sec:sim:routing}).
The postponement procedure of an order, see \Cref{alg:ord:opp}, is a simple method that yields an \texttt{Event} with the time until which the order has been postponed.
If the postponement procedure is processed, a callback on order postponement expiration is invoked, requesting routing by default.
Also recall that the postponement process will be interrupted if a routing is requested in the meantime (see Step~1 in \Cref{sec:sim:routing}).
These two cases correspond to Case~1 and Case~2 in \Cref{sec:bas:ord:post}, respectively.

\begin{algorithm}
\caption{Postponement procedure for orders (default)}\label{alg:ord:opp}
\begin{algorithmic}[1]
\State $t\gets$ the time until the order is postponed
\If{ $\tnow < t$ }
  \Yield $\operatorname{event}(t)$ \Comment{order postponement expired}
\EndIf
\end{algorithmic}
\end{algorithm}

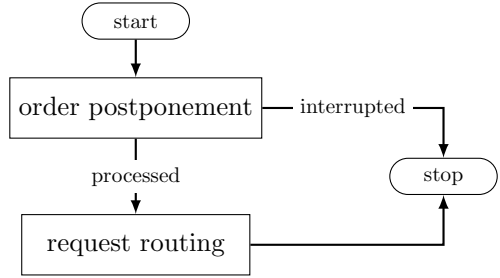
\begin{figure}
\centering
\begin{tikzpicture}
\node[fc-startstop] (start) at (0,0) {start};
\node[fc-state,below=0.5cm of start] (postponement) {order postponement};
\node[fc-state,below=1cm of postponement] (routing) {request routing};
\coordinate (cstop) at ($(postponement.east)!0.5!(routing.east)$);
\node[fc-startstop,right=1.75cm of cstop] (stop) {stop};

\draw[fc-arrow] (start) -- (postponement);
\draw[fc-arrow] (postponement) -- (routing) node[midway,fc-arrowlabel] {processed};
\draw[fc-arrow] (routing) -| (stop);

\draw[fc-arrow] (postponement.east) -| (stop.north) node[fc-arrowlabel,near start] {interrupted};

\end{tikzpicture}
\caption{
Postponement procedure for postponed orders.
}
\label{fig:sim:opp}
\end{figure}

\subsubsection{Order callbacks}

An \texttt{Order} has several callback methods that are invoked, for example, when the order is requested, rejected, canceled, postponed, picked up, delivered, or when the postponement of the order is expired.
By requesting routing in such a callback, we can model, for example, decision points on order request/cancellation/postponement.

\subsection{Vehicles}\label{sec:sim:vehicles}

\begin{figure}
\centering
\begin{tikzpicture}
\node[fc-startstop] (start) at (0,0) {start};
\node[fc-quest,below=0.5cm of start] (nextvisit) {vehicle has next visit};
\node[fc-state,below=1cm of nextvisit] (predeparture) {pre-departure};
\node[fc-state,below=1cm of predeparture] (travel) {travel};
\node[fc-state,below=1cm of travel] (preservice) {pre-service};
\node[fc-state,below=1cm of preservice] (service) {service};
\coordinate (cstop) at ($(preservice.east)!0.5!(travel.east)$);
\node[fc-startstop,right=1.75cm of cstop] (stop) {stop};

\draw[fc-arrow] (start) -- (nextvisit);
\draw[fc-arrow] (nextvisit) -- (predeparture) node[midway,fc-arrowlabel] {true};
\draw[fc-arrow] (predeparture) -- (travel) node[midway,fc-arrowlabel] {processed};
\draw[fc-arrow] (travel) -- (preservice) node[midway,fc-arrowlabel] {processed};
\draw[fc-arrow] (preservice) -- (service) node[midway,fc-arrowlabel] {processed};

\draw[fc-arrow] (service.west) -- ++(-1.75cm,0) node[fc-arrowlabel,midway] {processed} |- (nextvisit.west) ;

\draw[fc-arrow] (nextvisit.east) -| ($(stop.north)+(0.30,0)$) node[fc-arrowlabel,near start] {false};
\draw[fc-arrow] (predeparture.east) -| ($(stop.north)+(0.00,0)$) node[fc-arrowlabel,near start] {interrupted};
\draw[fc-arrow,dashed] (travel.east) -| ($(stop.north)-(0.30,0)$) node[fc-arrowlabel,near start] {interrupted};
\draw[fc-arrow,dashed] (preservice.east) -| ($(stop.south)-(0.15,0)$) node[fc-arrowlabel,near start] {interrupted};
\draw[fc-arrow,dashed] (service.east) -| ($(stop.south)+(0.15,0)$) node[fc-arrowlabel,near start] {interrupted};

\scoped[on background layer]
{
  \node[fc-fit,fit=(predeparture)(travel)(preservice)(service)] {};
}

\end{tikzpicture}
\caption{
Execution procedure of a vehicle.
Dashed arrows refer to non-default connections.
}
\label{fig:sim:vexe}
\end{figure}
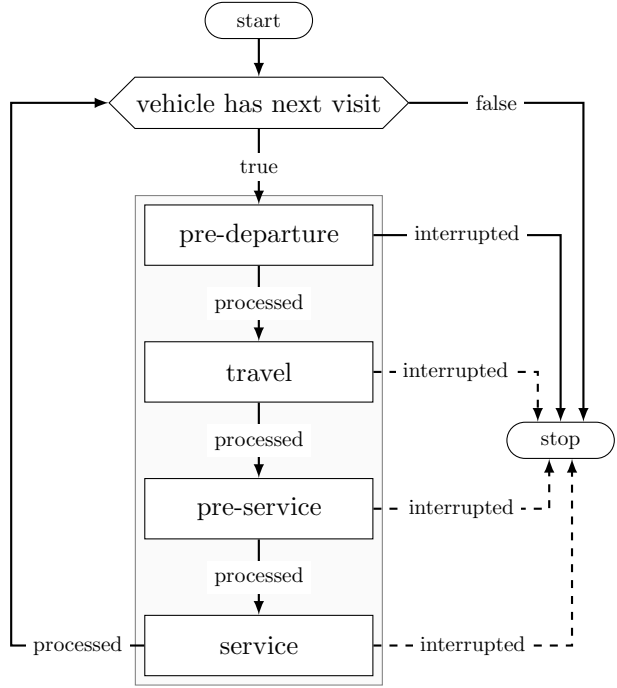

Each vehicle is modeled with an instance of the \texttt{Vehicle} class.
The \emph{execution procedure} of a vehicle executes its route plan as described in \Cref{sec:bas:veh:exe}.
This execution procedure consists of four consecutive parts, these are, the \emph{pre-departure procedure}, the \emph{departure procedure}, the \emph{pre-service procedure}, and the \emph{service procedure}, see \Cref{fig:sim:vexe}.
After the route plans are updated during a routing procedure (\Cref{sec:sim:routing}, Step~5), idle vehicles start their execution procedure (i.e., a procedure  for each idle vehicle is added to the simulation environment in a separate \texttt{Process}).

\subsubsection{Pre-departure procedure}\label{sec:sim:veh:predep}
The \emph{pre-departure procedure} models the phase before departure.
The \texttt{Event} indicating the end of the procedure corresponds to a vehicle departure event described in \Cref{sec:sdp:events}.
The procedure is suitable, for example, to delay the departure of the vehicles.

\begin{algorithm}
\caption{Pre-departure procedure (default)}\label{alg:veh:predep}
\begin{algorithmic}[1]
\State $t\gets$ earliest start time for the next visit (0, if not given)
\If{$\tnow < t$}
  \Yield $\operatorname{event}(t)$ \Comment{departure postponement expired}
  \State request routing (optional)
\EndIf
\Yield $\operatorname{timeout}(0)$ \Comment{vehicle departure}
\end{algorithmic}
\end{algorithm}

\paragraph{Delaying the departure}
The default pre-departure procedure, see \Cref{alg:veh:predep}, checks whether an earliest start time that has not yet passed is associated with the vehicle's next visit, and if so, yields an \texttt{Event} with that time.
By this, the travel procedure of the vehicle will start at that earliest start time.
Note, however, that the pre-departure procedure can be interrupted, as happens when a routing procedure is about to begin (\Cref{sec:sim:routing}, Step~1).
This behavior is consistent with the procedure proposed in the modeling framework (\Cref{sec:bas:veh:del}).

\subsubsection{Travel procedure}
The \emph{travel procedure} models the travel of the vehicle from its current location to its next location to visit.
By default, see \Cref{alg:veh:travel}, the procedure gets the corresponding travel time, and yields a timeout \texttt{Event} with that delay.
This \texttt{Event} corresponds to a vehicle arrival event described in \Cref{sec:sdp:events}.
The travel time procedure can be used, for example, to model deterministic, time-dependent, or stochastic travel times, etc.

\begin{algorithm}
\caption{Travel procedure (default)}\label{alg:veh:travel}
\begin{algorithmic}[1]
\State $\delta\gets$ travel time from current location to next location
\Yield $\operatorname{timeout}(\delta)$ \Comment{vehicle arrival}
\end{algorithmic}
\end{algorithm}

\paragraph{Deterministic and stochastic travel times}

The travel time between the corresponding locations can be a pre-calculated constant value or can be calculated on the fly.
Since the procedure has access to the current simulation time, it is easy to implement time-dependent travel times as well.
The procedure is also suitable for implementing stochastic travel times, see \Cref{alg:stoch:travel} for a primitive example.

\begin{algorithm}
\caption{Custom travel procedure with stochastic travel times}\label{alg:stoch:travel}
\begin{algorithmic}[1]
\State $\delta_1\gets$ travel time from current location to next location
\Yield $\operatorname{timeout}(\delta_1)$
\State $\delta_2\gets$ random delay
\Yield $\operatorname{timeout}(\delta_2)$ \Comment{vehicle arrival}
\end{algorithmic}
\end{algorithm}

\paragraph{En route diversion}
Our modeling framework does not allow en route diversion, so travel \texttt{Process}es cannot be interrupted by default.
However, all the necessary callbacks for this purpose are provided in our simulator.

\subsubsection{Pre-service procedure}
The \emph{pre-service procedure} models the phase between an arrival and the start of the subsequent service.
The \texttt{Event} indicating the end of the procedure corresponds to a service start event described in \Cref{sec:sdp:events}. 
The procedure can be used to delay the start of the service, and thus it is suitable, for example, for modeling earliest service start times, docking restrictions, stochastic completion times, etc.

\begin{algorithm}
\caption{Pre-service procedure (default)}\label{alg:veh:preserv}
\begin{algorithmic}[1]
\State $t\gets$ latest earliest pickup/delivery service start time of orders to pickup/deliver
\State $r\gets$ service request
\Yield $\operatorname{allof}(\operatorname{event}(t),r)$ \Comment{service start}
\end{algorithmic}
\end{algorithm}

\paragraph{Time windows}
Orders can be associated with \emph{earliest pickup/delivery service start times} (\Cref{sec:bas:orders}).
These are hard constraints, that is, even if the vehicle has arrived at the location, the service cannot start until the time window for each order tp deliver/pickup is opened.
The default pre-service procedure, see \Cref{alg:veh:preserv}, creates an \texttt{Event} that is processed at the latest earliest service start time, so sthat the service cannot start earlier.

\paragraph{Docking restrictions}
Shared resources can be associated with the locations to limit the number of vehicles that can be served simultaneously.
For example, \texttt{Resource}s can be used to model the number of docking ports at the locations.
The default pre-service procedure, see \Cref{alg:veh:preserv}, makes a request on the \texttt{Resource} of the corresponding location, so the service can only start after a docking port has become available for the vehicle.

\paragraph{Stochastic completion times}
In the restaurant meal delivery problem of \citep{ulmer2021restaurant}, the completion times of the orders are stochastic, and the drivers may need to wait at restaurants for the order to be ready.
The custom pre-service procedure depicted in~\Cref{alg:rmdp:preserv} models this behavior.

\begin{algorithm}
\caption{Custom pre-service procedure for stochastic ready times}\label{alg:rmdp:preserv}
\begin{algorithmic}[1]
\State $t\gets$ ready time of the order to pickup ($\tnow$, if it is ready)
\Yield $\operatorname{event}(t)$ \Comment{service start}
\end{algorithmic}
\end{algorithm}

\subsubsection{Service procedure}
The \emph{service procedure} models the service of the vehicle at its current location.
The main purpose of the procedure is to model the pickup and the delivery of the orders, see \Cref{alg:veh:service}, however, other vehicle activities (e.g., parking) can also be included.
The service procedure can be used, for example, to model deterministic service times, stochastic service times, order-independent service times, etc.

\begin{algorithm}
\caption{Service procedure (default)}\label{alg:veh:service}
\begin{algorithmic}[1]
\For{order~$o_i$ to deliver}
\State $\delta_d\gets$ delivery service time of the order \label{alg:veh:service:dstime}
\Yield $\operatorname{timeout}(\delta_d)$ \Comment{order delivery}
\EndFor
\For{order~$o_i$ to pickup}
\State $\delta_p\gets$ pickup service time of the order \label{alg:veh:service:pstime}
\Yield $\operatorname{timeout}(\delta_p)$ \Comment{order pickup}
\EndFor
\Yield $\operatorname{timeout}(0)$ \Comment{service finish}
\end{algorithmic}
\end{algorithm}

\paragraph{Deterministic and stochastic service times}
In the default service procedure, see \Cref{alg:veh:service}, unloading takes place first and then loading takes place afterwards.
Each pickup/delivery is modeled with a timeout \texttt{Event} with a delay corresponding to the duration of the pickup/delivery.
That service times can be either deterministic or stochastic, and can depend on the order, on the vehicle, and on the location as well.

\paragraph{Order-independent service times}
Order-independent service times (e.g., parking, docking) can be also taken into consideration.
For example, in case of the dynamic pickup and delivery problem of \citet{hao2022introduction}, the services start with a fixed-length docking, which can be modeled by the procedure depicted in \Cref{alg:dpdp:service}.

\begin{algorithm}
\caption{Custom service procedure}\label{alg:dpdp:service}
\begin{algorithmic}[1]
\State $\delta\gets$ dock approaching time
\Yield $\operatorname{timeout}(\delta)$
\State apply default service procedure
\end{algorithmic}
\end{algorithm}

\subsubsection{Callbacks}
\texttt{Vehicle}s have several callback methods, which are invoked, for example, when the vehicle arrives/departs at/from a location, when the service of the vehicle starts/finishes, or when one of its process is interrupted.
By requesting routing in such a callback, we can model, for example, decision points on vehicle arrival/departure, and decision points on service start/finish.

\end{document}